\RequirePackage{snapshot}
\documentclass[twoside]{article}
\usepackage[a4paper]{geometry}
\usepackage[latin1]{inputenc} 
\usepackage[T1]{fontenc} 
\usepackage{RR}
\usepackage{hyperref}
\RRNo{8602}
\RRdate{September 2014}
\RRauthor{
Vincent Acary\thanks{vincent.acary@inria.fr}}
\authorhead{Acary}
\RRtitle{Conservation d'énergie et propriétés de dissipation de méthodes d'intégration en temps pour la dynamique élastique non--régulière avec contact}
\RRetitle{Energy conservation and dissipation properties of time-integration methods\\ for the nonsmooth elastodynamics with contact}
\titlehead{Energy conservation and dissipation properties for the nonsmooth elastodynamics with contact}
\RRresume{Ce rapport de recherche propose une étude des propriétés de conservation et de dissipation de l'énergie mécanique pour différents schémas d'intégration en temps de la dynamique élastique avec du contact unilatéral. Sachant que l'application directe des schémas standards de type Newmark et des schémas $\alpha$--généralisés conduisent à des explosions de l'énergie mécanique, on étudie deux schémas dédiés à l'intégration en temps des systèmes non réguliers avec contact : le schéma de Moreau--Jean et le schéma $\alpha$--généralisé non--régulier. La conservation de l'énergie et les propriétés de dissipation du schéma de Moreau--Jean sont d'abord démontrées. Dans un second temps, le schéma $\alpha$--généralisé non--régulier est étudié en adaptant les travaux précurseurs de Krenk et  H{\o}gsberg dans le contexte du contact unilatéral. Finalement, les propriétés connues du schéma de Newmark et du schéma Hilber--Hughes--Taylor (HHT) dans le cas régulier sont étendues dans le cas avec contact sans hypothèses supplémentaires.  }
\RRabstract{
This research report is devoted to the study of the conservation and the dissipation properties of the mechanical energy of several time--integration methods dedicated to the elasto--dynamics with unilateral contact.  Given that the direct application of the standard schemes as the Newmark schemes or the generalized--$\alpha$ schemes leads to energy blow-up, we study two schemes dedicated to the time--integration of nonsmooth systems with contact:  the Moreau--Jean scheme and the nonsmooth generalized--$\alpha$ scheme.  The energy conservation and dissipation properties of the Moreau--Jean is firstly shown.  In a second step,  the nonsmooth generalized--$\alpha$ scheme is studied by adapting the previous works of Krenk and  H{\o}gsberg in the context of unilateral contact. Finally, the known properties of the Newmark  and the Hilber--Hughes--Taylor (HHT) scheme  in the unconstrained case are extended without any further assumptions to the case with contact.  }
\RRmotcle{Dynamique du contact, impact, mécanique numérique du contact, intégration numérique en temps, conservation d'énergie, propriétés de dissipation.}
\RRkeyword{contact dynamics, impact,  computational contact mechanics,  numerical time integration, energy conservation, dissipation properties}
\RRprojet{Bipop}  
 \RCGrenoble 

\usepackage{epsfig,amsmath,amsfonts,a4wide,cases,amsmath,amssymb,subfigure,color,hyperref,algorithm,algorithmic,color,natbib,url}
\usepackage[squaren]{SIunits}

\newcommand{\Frac}[2]{{\displaystyle \frac{\displaystyle #1}{\displaystyle #2}}}

\def\XXint#1#2#3{{\setbox0=\hbox{$#1{#2#3}{\int}$}
     \vcenter{\hbox{$#2#3$}}\kern-.5\wd0}}

\def\nbC{{\mathchoice {\setbox0=\hbox{$\displaystyle\rm C$}%
\hbox{\hbox to0pt{\kern0.4\wd0\vrule height0.9\ht0\hss}\box0}}
{\setbox0=\hbox{$\textstyle\rm C$}\hbox{\hbox
to0pt{\kern0.4\wd0\vrule height0.9\ht0\hss}\box0}}
{\setbox0=\hbox{$\scriptstyle\rm C$}\hbox{\hbox
to0pt{\kern0.4\wd0\vrule height0.9\ht0\hss}\box0}}
{\setbox0=\hbox{$\scriptscriptstyle\rm C$}\hbox{\hbox
to0pt{\kern0.4\wd0\vrule height0.9\ht0\hss}\box0}}}}
\def\nbQ{{\mathchoice {\setbox0=\hbox{$\displaystyle\rm
Q$}\hbox{\raise
0.15\ht0\hbox to0pt{\kern0.4\wd0\vrule height0.8\ht0\hss}\box0}}
{\setbox0=\hbox{$\textstyle\rm Q$}\hbox{\raise
0.15\ht0\hbox to0pt{\kern0.4\wd0\vrule height0.8\ht0\hss}\box0}}
{\setbox0=\hbox{$\scriptstyle\rm Q$}\hbox{\raise
0.15\ht0\hbox to0pt{\kern0.4\wd0\vrule height0.7\ht0\hss}\box0}}
{\setbox0=\hbox{$\scriptscriptstyle\rm Q$}\hbox{\raise
0.15\ht0\hbox to0pt{\kern0.4\wd0\vrule height0.7\ht0\hss}\box0}}}}
\def\nbT{{\mathchoice {\setbox0=\hbox{$\displaystyle\rm
T$}\hbox{\hbox to0pt{\kern0.3\wd0\vrule height0.9\ht0\hss}\box0}}
{\setbox0=\hbox{$\textstyle\rm T$}\hbox{\hbox
to0pt{\kern0.3\wd0\vrule height0.9\ht0\hss}\box0}}
{\setbox0=\hbox{$\scriptstyle\rm T$}\hbox{\hbox
to0pt{\kern0.3\wd0\vrule height0.9\ht0\hss}\box0}}
{\setbox0=\hbox{$\scriptscriptstyle\rm T$}\hbox{\hbox
to0pt{\kern0.3\wd0\vrule height0.9\ht0\hss}\box0}}}}
\def\nbS{{\mathchoice
{\setbox0=\hbox{$\displaystyle     \rm S$}\hbox{\raise0.5\ht0%
\hbox to0pt{\kern0.35\wd0\vrule height0.45\ht0\hss}\hbox
to0pt{\kern0.55\wd0\vrule height0.5\ht0\hss}\box0}}
{\setbox0=\hbox{$\textstyle        \rm S$}\hbox{\raise0.5\ht0%
\hbox to0pt{\kern0.35\wd0\vrule height0.45\ht0\hss}\hbox
to0pt{\kern0.55\wd0\vrule height0.5\ht0\hss}\box0}}
{\setbox0=\hbox{$\scriptstyle      \rm S$}\hbox{\raise0.5\ht0%
\hboxto0pt{\kern0.35\wd0\vrule height0.45\ht0\hss}\raise0.05\ht0%
\hbox to0pt{\kern0.5\wd0\vrule height0.45\ht0\hss}\box0}}
{\setbox0=\hbox{$\scriptscriptstyle\rm S$}\hbox{\raise0.5\ht0%
\hboxto0pt{\kern0.4\wd0\vrule height0.45\ht0\hss}\raise0.05\ht0%
\hbox to0pt{\kern0.55\wd0\vrule height0.45\ht0\hss}\box0}}}}
\def\nbZ{{\mathchoice {\hbox{$\sf\textstyle Z\kern-0.4em Z$}}
{\hbox{$\sf\textstyle Z\kern-0.4em Z$}}
{\hbox{$\sf\scriptstyle Z\kern-0.3em Z$}}
{\hbox{$\sf\scriptscriptstyle Z\kern-0.2em Z$}}}}

\def\n{{\hbox{\tiny{N}}}}
\def\t{{\hbox{\tiny{T}}}}
\def\ss{{\hbox{\tiny{S}}}}

\def\m{\,\mathsf{m}}

\def\N{\,\mathsf{N}}

\def\p{\partial}

\def\div{\mathop{\rm div}\nolimits}
\def\ae#1{\stackrel{\mbox{\scriptsize a.e.}}{#1}}

\def\eqref#1{{\rm (\ref{#1})\/}}
\def\p{\mathord{\rm proj}}
\def\N{\mathord{\rm N}}

\def\qedtext{\mbox{}\hfill$\Box$}

\def\real{\mbox{\scriptsize I\hspace{-.15em}R}}

\catcode`\é=13
\defé{\'e}
\catcode`\è=13
\defè{\`e}
\catcode`\à=13
\defà{\`a}
\catcode`\ç=13
\defç{\c c}
\def\N{\mbox{I\hspace{ -.15em}N}}

\def\LMGC90{{\small \it LMGC90 }}

\def\12T{{\frac{_{1}}{^{2T}}}}

\def\geq{\geqslant}
\def\leq{\leqslant}

\begingroup
\def\2#1{\ifnum#1<10 0\fi\the#1}
\xdef\isodayandtime{\the\year-\2\month-\2\day\space\2{\count0}:%
\2{\count2}}
\endgroup

\makeatletter

{\newtheorem{ndr1bb}{\textbf{\textsc{Redaction note B.B.}}}[section]}

{%
\noindent\begin{ndr1bb}\hrule\vspace{1em}%
\ttfamily\small
}%
{%
\begin{flushright}%
\end{flushright}%
\vspace{-1.5em}\hrule
\end{ndr1bb}%
}

{\newtheorem{ndr1va}{\textbf{\textsc{Redaction note V.A.}}}[section]}

{%
\noindent\begin{ndr1va}\hrule\vspace{1em}%
\ttfamily\small \  \\
\indent}%
{%
\begin{flushright}%
\  \\
\end{flushright}%
\vspace{-1.5em}\hrule
\end{ndr1va}%
}
\makeatother

 \def\c{\mathop{{\rm 1}\mskip-10.0mu{\rm C}}\nolimits}

\newcommand{\RR}{\mbox{\rm $I\!\!R$}}

\newtheorem{proposition}{Proposition}
\newtheorem{lemma}{Lemma}

\newtheorem{remark}{Remark}

\newtheorem{corollary}{Corollary}

\def\dt{{\rm d}t}
\def\dv{{\rm d}v}
\def\di{{\rm d}i}
\def\dI{{\rm d}I}

\def\newW{{W}}

\begin{document}
\makeRR   
\tableofcontents

\section{Introduction and motivations}

The numerical time integration of mechanical systems with unilateral contact is known to be a difficult task, mainly due to the nonsmoothness of the dynamic response when a contact occurs. For two recent reviews of the existing methods in the literature and the associated issues, we refer to~\citep{Doyen.ea2011,Krause.Walloth2012} and the standard textbooks~\citep{Laursen2003,Wriggers2006}. One of the main conclusions  is that standard schemes widely used in computational contact Mechanics, such as the Newmark scheme, the Hilber--Hughes--Taylor scheme (HHT) or the generalized--$\alpha$ scheme cannot be directly applied to the simulation of systems with unilateral contact and impact.

 In the most favorable cases, these schemes exhibit artificial oscillations of the contact velocities and forces, that blur the whole stresses in the structure. The source of these artificial oscillations is the nonsmoothness of the contact conditions that yields a jump in the velocities when a contact is closing. In finite--freedom Mechanics, when one deals with space--discretized structures after a semi--discretization, the velocity jump of a finite mass is associated with an impact. Hence, an impact law has to be specified to close the equations of the system. The low regularity of the velocity and the distributional character of the reaction prevents the use of schemes with an high order of accuracy. A partial remedy for theses problems is to use a fully implicit first--order approximation of the contact forces as it was suggested in~\citep{Jean.Moreau1987,Moreau1988,Carpenter.Taylor.ea92,Kane.ea1999} and/or a treatment of the constraints at the velocity level together with an (possibly perfectly inelastic or implicitly defined) impact law~\citep{Moreau1988,Laursen.Chawla97}. This work has been extended to nonlinear elastodynamics by~\citep{Hauret.LeTallec2006}. In~\citep{Ayyad.ea2009}, the nonlinear elastodynamics with a constraint at the position level is also considered. This yields a similar numerical scheme for the energetic properties but a additional step, in the same vein as in \citep{Laursen.Love2002}, is added to correct the constraint at the position level.


 In the worst cases, the standard schemes exhibit numerical instabilities, and possibly, energy blow-ups. An alternative solution is to design numerical schemes that conserve or dissipate energy. In~\citep{Laursen.Chawla97}, the energy--conserving scheme based on the pioneering works Simo and co-workers~\citep{Simo.Tarnow92,Simo.ea1995} for nonlinear elastodynamics is extended to the elasto--dynamical with unilateral contact. This work results in the use of the mid--point rule together with a velocity--level constraints and an implicit treatment of the contact forces.

The velocity--level formulation has two main advantages: it controls the dissipation of energy  at contacts and it stabilizes the contact velocity. One of the main drawback is the violation of the constraints on the position level which is proportional to the time--step. In \citep{Laursen.Love2002}, the authors propose to satisfy the constraints at the position level together the energy conservation by introducing an artificial velocity variable at the price to have oscillations of the contact velocity. In the latter case, the velocity--level constraints is not satisfied. Alternatively, the constraints at the velocity level and at the position level can be both satisfied by adding an artificial multiplier to perform a projection on the position constraint~\citep{Acary2013224}. As remarked in~\citep{Krause.Walloth2012}, there exists no algorithm satisfying the constraints at both position and velocity levels and ensuring the energy conservation simultaneously.

In the previous attempts to adapt the Newmark--based schemes in the context of computational contact mechanics,  very few results are available on the energy conservation or dissipation when we deal with the  unilateral contacts at the position level. In the unconstrained case, the classical Newmark scheme, with the special choice of parameters $\gamma=2\beta=1/2$, leads to an algorithm conserving the total energy of the system. With  position based unilateral constraints, it is shown in \citep{Krause.Walloth2012} that the scheme cannot conserve the energy even if $\gamma=2\beta=1/2$. In~\citep{Khenous2005}, the author can only conclude to the energy dissipation of the Newmark scheme is the special case of a very dissipative order one scheme with $\gamma = 2\beta=1$. With a full implicit treatment of the constraint as in \citep{Kane.ea1999}, as far as we know there is no general study of the energy properties for all admissible values of $\gamma$ and $\beta$. When we consider the Newmark scheme with a velocity--level formulation of the constraints, the only available results  are those in~\citep{Laursen.Chawla97} that can be adapted to the special case of the Newmark scheme with a fully implicit treatment of the Lagrange multiplier. Indeed, the midpoint rule together with a velocity--level formulation is very similar to the Moreau--Jean scheme~\citep{Jean.Moreau1987,Moreau1988,Moreau1999,Jean1999} based on the $\theta-$method when $\theta$ is equal to $1/2$. This latter scheme is in turn equivalent to the Newmark scheme with $\gamma=2\beta=1/2$.  Besides these special cases, the general case of the HHT scheme and generalized--$\alpha$ scheme  are  not treated from the energy properties point of view. This report attempts to bridge this gap.

Very recently, a new class of schemes has been proposed in~\citep{Chen.ea2012,Chen.ea_IJNME2013} which takes advantage of the Moreau--Jean scheme in terms of robustness and stability while adding some key properties of the Newmark--based schemes, that are the second--order approximation of the smooth force terms and the controlled damping of the high--frequency dynamics. This work yields the so--called nonsmooth Newmark, nonsmooth HHT and nonsmooth generalized--$\alpha$ schemes that deal with the contact forces trough their associated impulses in a fully implicit way, and treat the constraints at the velocity level together with Newton's impact law.  In this report, the main goal is to show that the Moreau--Jean scheme and the nonsmooth schemes have the same energetic behavior as their counterparts in the unconstrained case. To this aim, the detailed list of the objectives is as follows:
\begin{itemize}
\item to show the energy conserving and decaying properties of the Moreau--Jean scheme.
\item to provide results on the algorithmic energy conservation and dissipation of the nonsmooth Newmark scheme.
\item to apply and to extend the techniques developed in~\citep{Krenk.Hogsberg2005,Krenk2006} to study the dissipation properties of the nonsmooth generalized--$\alpha$ schemes. 
\item to show by means of the previous method that the nonsmooth HHT scheme dissipates a kind of algorithmic energy.
\item to propose an alternative $\alpha$--scheme as in~\citep{Krenk.Hogsberg2005,Krenk2006} that dissipates a kind of  algorithmic energy.
\end{itemize}

The report is organized as follows. Section~\ref{Sec:Nonsmoothmodelling} recalls the basic ingredients of the nonsmooth modeling of finite--dimensional mechanical systems subjected to unilateral contact and impact. In Section~\ref{Sec:EnergyBalance}, the energy balance of a mechanical system with jumps in the velocity is formulated. The main schemes, studied in this report, are detailed in Section~\ref{Sec:Background}. The energy analysis of the Moreau--Jean scheme is done in Section~\ref{Sec:EnergyAnalysisMoreau}. Section~\ref{Sec:Dissipation-Alpha-Scheme} starts by the presentation of the Krenk--H{\o}gsberg method for the analysis of the discrete energy balance over a time--step for the $\alpha$--schemes. After a first general result on the nonsmooth generalized--$\alpha$ schemes and the alternative nonsmooth Krenk--H{\o}sberg  generalized--$\alpha$ scheme, the nonsmooth HHT case  and the nonsmooth Newmark case are fully developed. Section~\ref{Sec:Conclusion} concludes the report.

\paragraph{Notation}  The following notation is used throughout the paper. The Euclidean norm for a vector $x\in\RR^n$ is denoted by $\|x\|$. For positive definite (respectively positive semi--definite) matrix $M\in \RR^{n\times n}$,  $\|x\|_M$ denotes the norm (respectively the semi--norm) in the metric defined by $M$.
Let $I$ denote a real time interval of any sort. For a function $f : I \rightarrow \RR^n$ of Bounded Variation (BV), we denote the right--limit function by $f^+(t) = \lim_{s\rightarrow t, s>t} f(s)$, and respectively  the left--limit by $f^-(t) = \lim_{s\rightarrow t, s<t} f(s) $. 
 We denote by $0= t_0 <t _1 <\ldots < t_k  < \ldots < t_N =T $ a finite partition (or a subdivision) of the  time interval $ [0,T] \; (T > 0)$. The integer $N$ stands for the number of time intervals in the subdivision.  The length of a time step is denoted by $h_k = t_{k+1}- t_k$.  For simplicity's sake, the schemes are  presented in the sequel with a time step  denoted by $h$ for short. The  value of a real function $x(t)$  at the time $t_k$, is approximated by $x_k$. In the same way, the  notation $x_{k+\theta} = (1-\theta)x_{k}+\theta x_{k+1}$ is used for $\theta \in [0,1]$. The following notation is introduced to analyze the energetic behavior of the nonsmooth generalized-$\alpha$ scheme
\begin{equation}
  \label{eq:notation}
  \begin{array}{lcl}
  x_{k,\gamma} = \gamma x_{k+1} + (1-\gamma) x_{k}, & &
  x_{k-1,\gamma} = \gamma x_{k} + (1-\gamma) x_{k-1}.\\
\end{array}
\end{equation}
This notation generalizes the previous notation $x_{k+\theta} = (1-\theta)x_{k}+\theta x_{k+1}$ to avoid the ambiguity when a multi--step method is studied. For a function $f: \RR \rightarrow \RR^n$ and $h>0$,  one writes $f(h) = \mathcal O(h)$ if and only if there exist positive numbers $\delta$ and $M$ such that $\|f(h)\| \leq M h $ for $h< \delta$.
The notation $\dt$ defines the Lebesgue measure on $\RR$.

\section{Nonsmooth mechanical systems with unilateral contact}
\label{Sec:Nonsmoothmodelling}

Let us consider the equations of motion of a mechanical system subjected to unilateral constraints in the linear case:
\begin{subnumcases}{
    \label{eq:Lagrange-uni}}
  q(t_0) = q_0,\, v(t_0)=v_0, \\[1mm]
  \dot q(t) = v(t), \label{eq:Lagrange-uni-velo} \\[1mm]
  M \dot v(t)  + K q(t) + C v(t) =   G \lambda(t), \\[1mm]
  g(q(t)) = G^\top q(t) + w \geq 0, \quad \lambda(t) \geq 0,\quad  g^\top(q(t)) \lambda(t) =0,\label{eq:Lagrange-uni-d}
 \end{subnumcases}
where
\begin{itemize}
\item $q(t)\in \RR^n$ is the vector of  generalized coordinates and $v(t) = \dot q(t)$ the associated vector of generalized velocities,
\item the initial conditions are $q_0\in \RR^n$ and $v_0\in \RR^n$,
\item $M \in \RR^{n\times n}$ is the symmetric inertia matrix,  $K\in \RR^{n\times n}$ and $C\in \RR^{n\times n}$ are respectively  the stiffness  and the damping matrices,
\item the function $g(q(t))\in \RR^m$, called the gap function is used to define the unilateral constraints; with a slight abuse of notation we will also write the gap function as $g(t)$,
\item  the Jacobian matrix of $g$ with respect to $q$ is  $G^\top = \nabla^\top_q g(q(t))$ and is assumed to be constant in the linear setting, $w\in \RR^m$ is a constant vector,
\item $\lambda(t) \in \RR^m$ is the Lagrange multiplier vector associated with the constraints.
\end{itemize}

The constitutive law~(\ref{eq:Lagrange-uni-d}) for the perfect unilateral constraints is also termed the Signorini condition and can be written equivalently as
 \begin{equation}
   0 \leq g(t) \perp \lambda(t) \geq 0,
   \label{eq:signo}
 \end{equation}
where the inequalities involving vectors are understood to hold component-wise and the $\perp$ symbol means that $y^\top \lambda =0$. Finally, let us define the following variables relative to the constraints, called \textit{local variables}: the local velocity $U(t)$  and the (local) Lagrange multiplier $\lambda(t)$ which is associated with the generalized reaction forces $r(t)$ such that 
\begin{equation}
  \label{eq:lambda-gene}
  U(t) = \dot g(q(t))=  G^\top \, v(t), \quad\quad  r(t) = G \lambda(t).
\end{equation}

For finite-freedom mechanical systems with unilateral constraints as in  (\ref{eq:Lagrange-uni}), it is well--known that impacts may occur if the relative velocity at contact $U(t)$ is not compatible with the constraints. In other terms, if the contact indexed by $\alpha$ is closing with a negative relative velocity at $t_\star$, that is $U^{\alpha,-}(t_\star)<0$, the velocity has to jump in order to satisfy the constraints after the time of impact. However, the right velocity at the impact $U^{\alpha,+}(t_\star)$ is not determined by the system~(\ref{eq:Lagrange-uni-d}). This is the reason why 
 an impact law must be added to close the system of equations. In particle or rigid body dynamics,  the most simple impact law  is Newton's impact law
\begin{equation}
  \label{eq:newtonlaw}
   U^{\alpha,+}(t) = - e^\alpha \,U^{\alpha,-}(t), \text{ if } g^\alpha(t) =0 \text{ and } U^{\alpha,-}(t) \leq 0 , \quad \alpha \in \mathcal I
\end{equation}
where $e^\alpha$ is the coefficient of restitution at contact $\alpha$ and  $\mathcal I =\{1\ldots m\}$ is  the set of indices of constraints.

Usually, the condition in~(\ref{eq:newtonlaw})  only involves the condition $g^\alpha(t) =0$. Indeed, if the condition  $g^\alpha(t) \geq 0$ is always satisfied on the time interval of study, the relative pre-impact velocity $U^{\alpha,-}(t)$  cannot be positive when $g^\alpha(t) =0 $, except at the initial time. We will see further that this is not the case when the condition~(\ref{eq:newtonlaw}) is only prescribed at discrete time instants. We will have therefore to adapt the condition~(\ref{eq:newtonlaw}) in a suitable way in Section~\ref{Sec:Background} when  a time--discretization is performed.

The fact that the velocity $v(t)$ may encounter jumps yields some difficulties to define the acceleration everywhere in time. It is usual to assume that the velocity is a function of bounded variation that admits a right and left limit everywhere and which can be associated with a differential measure $\dv$ (see~\citep{Moreau1988b} for details). Almost everywhere with respect to the Lebesgue measure $\dt$, the velocity is differentiable in the classical way with respect to time and we get $\Frac{\dv}{\dt} = \ddot q(t)$. When a jump occurs, the standard differentiation cannot be applied since the acceleration is given by a Dirac distribution. Let us write the equation of motion in terms of differential measure in the linear case
\begin{eqnarray}
  \label{eq:signo-velocity}
  \begin{cases}
    M \dv  + K q(t)\,\dt + C v(t)\, \dt = F(t)\,\dt + G \,\dI, \\
    \dot q(t) = v^+(t), \\
    U(t) = G^\top v(t), \\
    g(t) = G^\top q(t) + w, \\
    \text{if }  g^\alpha(t) \leq 0 \text{ and } U^-(t) <0  , \text{ then } 0 \leq U^+(t) + e U^-(t) \perp \,\dI \geq 0,
  \end{cases}
\end{eqnarray}
where  $\di$ is the impulse reaction measure. Its local variant is defined by $\di = G \,\dI$. When the evolution is smooth, $r(t)$ is considered as the density of $\di$ with respect to the Lebesgue measure, that is
\begin{equation}
  \label{eq:density}
  r(t)  = \Frac{\,\di}{\dt}(t),\text{ or equivalently } \lambda(t)  = \Frac{\,\dI}{\dt}(t) .
\end{equation}
The last line of~(\ref{eq:signo-velocity}) defines the second--order Moreau sweeping process~\citep{Moreau1988}. It can be interpreted as a reformulation of the unilateral constraint~(\ref{eq:signo}) at the velocity level together with an impact law. It can also be equivalently viewed as an index--reduction technique, standard in the Differential Algebraic Equations (DAE) Theory when a constraint is differentiated to decrease the index of the DAE. It also contains to the so-called persistency conditions, defined in~\citep{Laursen.Chawla97} when the contact is closed, but adds the Newton impact law when a new contact occurs.

It can also be shown that the system~(\ref{eq:signo-velocity}) contains in a single formulation the nonimpulsive dynamics and the impact dynamics. If we omit the singular part in the decomposition of the measure, we can split the measure $\dv$ and $\di$ as follows
  splitting of measures
  \begin{equation}
    \label{eq:split}
    \begin{array}{lcl}
    \dv &=& \ddot q(t) \,\,\dt + \sum _i (v^+(t_i) -v^-(t_i) ) \delta_{t_i},\\ 
    \dI &=& \lambda \,\,\dt + \sum _i P_i\delta_{t_i}.
  \end{array}
  \end{equation}
Inserting this decomposition in the measure differential system~(\ref{eq:signo-velocity}), we get the standard non--impulsive equation of motion~(\ref{eq:Lagrange-uni}) almost everywhere and the impact equations at the time of impacts:
\begin{eqnarray}
  \label{eq:impact-velocity}
  \begin{cases}
    M  (v^+(t_i) -v^-(t_i) )   = G \,P_i, \\
    U^\pm(t_i) = G^\top v^\pm(t_i), \\
    g(t) = G^\top q(t) + w, \\
    \text{if } U^-(t_i) <0  , \text{ then } 0 \leq U^+(t_i) + e U^-(t_i) \perp P_i \geq 0.
  \end{cases}
\end{eqnarray}

 For more details on the modeling of multibody systems with unilateral constraints, we refer to~\citep{Acary.Brogliato2008bis,Pfeiffer.Glocker1996,Moreau1988} and for the mathematical analysis, we refer to \citep{Schatzman1978,Marques1993,Stewart2000a,Ballard2000}. 

\section{Energy balance analysis}
\label{Sec:EnergyBalance}
In the case on nonsmooth motion with impact, we recall that the equation of motion in terms of measures are given by
\begin{equation}
  \label{eq:LTI}
  \begin{cases}
    M \,\dv  + (K q(t) + C v(t))\,\dt = F \,\,\dt + \,\di, \\
    dq = v(t) \dt. 
  \end{cases}
\end{equation}
A detailed analysis of the energy balance for nonsmooth systems can be found in~\citep{BB-RL-BM-OE2007} and \citep{Leine.vandeWouw2008}. Let us recall in this paper the basic formulae. The energy balance can be obtained by multiplying the equation of motion by $v^++v^-$
\begin{equation}
  \label{eq:LTI-Energy-balance}
  \begin{array}{rcl}
   (v^+ + v^-)^\top M \,\dv  +  (v^+ + v^-)^\top  (Kq + C v)\,\dt &=&  (v^+ + v^-)^\top F \,\,\dt + (v^+ + v^-)^\top \,\di,
 \end{array}
\end{equation}
that is equivalently
\begin{equation}
  \label{eq:LTI-Energy-balance0}
  \begin{array}{rcl}
   \,{\rm d}(  v^\top M v)  +  (v^+ + v^-)  (Kq + C v)\,\dt &=&  (v^+ + v^-) F \,\,\dt + (v^+ + v^-) \,\di.
  \end{array}
\end{equation}
With the standard definition of the total mechanical energy of the system,
\begin{equation}
  \label{eq:LTI-energy-balance2}
  \begin{array}{rcl}
    \mathcal E  :=   \Frac 1 2 v^\top M v  +  \Frac 1 2 q^\top K q,
  \end{array}
\end{equation}
 we get the following energy balance
\begin{equation}
  \label{eq:LTI-energy-balance1}
  \begin{array}{rcl}
 2 {\rm d}{\mathcal E}  :=  \,d(  v^\top M v)  +  2 q^\top K {\rm d}q  &=&  2 v^\top F \,\,\dt  - 2v^\top C v \,\,\dt +  (v^+ + v^-)^\top \,\di.  
  \end{array}
\end{equation}
If we split the differential measure in $\di =  r(t) \,\,\dt  + \sum _i p_i \delta_{t_i}$, we get 
\begin{equation}
  \label{eq:LTI-energy-balance10}
  \begin{array}{rcl}
 2 {\rm d}{\mathcal E}  :=  \,{\rm d}(  v^\top M v)  +  2 q^\top K {\rm d}q  &=&  2 v^\top (F +r )\,\,\dt  - 2v^\top C v \,\,\dt +  \displaystyle\sum _i (v^+ + v^-)^\top p_i\delta_{t_i}.
  \end{array}
\end{equation}
By integration over a time interval $[t_0,t_0]$ such that $t_i \in (t_0,t_1)$, we obtain an energy balance equation (also named the dissipation equality~\citep{BB-RL-BM-OE2007}) as 
\begin{equation}
  \label{eq:LTI-energy-balance11}
  \begin{array}{rcl}
    \Delta \mathcal E =   \displaystyle \mathcal E(t_1) -\mathcal E(t_0) &=&     \displaystyle \int_{t_0}^{t_1} v^\top (F + r )\,\,\dt  \,-  \int_{t_0}^{t_1} v^\top C v \,\,\dt + \sum_i \Frac 1 2 (v^+(t_i) + v^-(t_i))^\top p_i.
  \end{array}
\end{equation}
The right hand side of the energy balance equation represents the work  done in the time interval $[t_0,t_1]$ that can be decomposed as follows:
\begin{itemize}
\item the term
  \begin{equation}
    \label{eq:LTI-energy-balance5}
   W^{\sf ext}=  \displaystyle \int_{t_0}^{t_1} v^\top F \,\,\dt,
  \end{equation}
is the work done by the external forces,
\item the term
  \begin{equation}
    \label{eq:LTI-energy-balance6}
   W^{{\sf damping}}=  - \displaystyle \int_{t_0}^{t_1} v^\top C v \,\,\dt,
  \end{equation}
is the work done by the damping internal forces. If the damping matrix $C$ is a positive semi--definite matrix, we then conclude on the sign of the work, $W^{{\sf damping}} \leq 0$.
\item the term
  \begin{equation}
    \label{eq:LTI-energy-balance7}
   W^{\sf contact}=  \displaystyle \int_{t_0}^{t_1} v^\top r \,\,\dt,
  \end{equation}
is the work done by the contact forces. If we consider perfect unilateral constraints, we have
  \begin{equation}
    \label{eq:LTI-energy-balance8}
    v^\top r = v^\top G \lambda = U^\top \lambda =0 
  \end{equation}
and then $W^{\sf contact}=0$,
\item the term 
  \begin{equation}
    \label{eq:LTI-energy-balance9}
    W^{\sf impact} = \Frac 1 2 (v^+(t_i) + v^-(t_i))^\top p,
  \end{equation}
represents the
  work done by the contact impulse at the time of impact $t_i$. Since $p =G P $  and $U=G^\top v $ , using the
  Newton impact law, we get
  \begin{equation}
    \label{eq:LTI-energy-balance3}
    \begin{array}{rcl}
     W^{\sf impact}=  \Frac 1 2 (v^+(t_i) + v^-(t_i))^\top p 
      &=&   \Frac 1 2 (v^+(t_i) + v^-(t_i))^\top G P\\[2mm]
      &=&   \Frac 1 2 (U^+(t_i) + U^-(t_i))^\top P\\[2mm]
      &=&   \Frac 1 2 ((1-e) U^-(t_i))^\top P \leq 0 \text { for   } 0\leq e \leq 1.
    \end{array}
  \end{equation}
  The formula~(\ref{eq:LTI-energy-balance9}) of the work done by the impulse during an impact is also known as Thomson and Tait's formula~\citep[Section 4.2.12]{Brogliato1999}. The fact that the work has to be negative such that an impact dissipates some energy is also related the Clausius--Duhem inequality applied to an impulse motion~\citep{Fremond2002}.
\end{itemize}



\section{Background on the time--integration methods}
\label{Sec:Background}
Leaving aside the time--integration methods based on an accurate event detection procedure (event--tracking schemes or event--driven schemes \cite[Chap. 8]{Acary.Brogliato2008bis}), the time--integration of nonsmooth mechanical systems is performed by means of event-capturing time--stepping schemes. In these schemes, the impact are not accurately located but captured by the refinement of the time-step size. One of the most well-proven method is the Moreau--Jean scheme~\citep{Jean.Moreau1987,Moreau1988,Jean1999}. This method which is sound from the mathematical analysis point of view (convergence analysis can be found in~\citep{Marques1993,dzonou2007,Dzonou.ea2009}) and which takes advantage of a strong practical experience, is of low order of accuracy, but very robust in many practical applications. Very recently, an attempt has been made in~\citep{Chen.ea2012,Chen.ea_IJNME2013} to improve its accuracy at least on the smooth terms of the equations of motion. This work has led to the nonsmooth Newmark and generalized-$\alpha$ scheme as it extends the standard schemes for computational mechanics to the nonsmooth dynamical case. Both schemes will be briefly presented in the following sections.

\subsection{Moreau--Jean's scheme}
 The Moreau--Jean scheme~\citep{Jean.Moreau1987,Moreau1988,Jean1999} is based on the Moreau sweeping process~(\ref{eq:signo-velocity}). The numerical time integration of~(\ref{eq:signo-velocity}) is performed on an interval $(t_k,t_{k+1}]$ of length $h$ as follows
 \begin{subnumcases}{\label{eq:MoreauTS}}
   M (v_{k+1}-v_k)
   + h K q_{k+\theta} + h C v_{k+\theta}     -   h F_{k+\theta} = p_{k+1} = G P_{k+1},\label{eq:MoreauTS-motion}\\[1mm] 
   q_{k+1} = q_{k} + h v_{k+\theta}, \quad \\[1mm]
   U_{k+1} = G^\top\, v_{k+1}, \\[1mm]
   \begin{array}{lcl}
     0 \leq U^\alpha_{k+1} + e  U^\alpha_{k} \perp P^\alpha_{k+1}  \geq 0,& \quad&\alpha \in \mathcal I_1, \\[1mm]
     P^\alpha_{k+1}  =0,&\quad& \alpha \in \mathcal I \setminus \mathcal I_1,
   \end{array}
 \end{subnumcases} 
with  $\theta \in [0,1]$. The following approximations are considered:
 \begin{eqnarray}
   \label{eq:no}
   v_{k+1}\approx v^+(t_{k+1});  \quad  U_{k+1}\approx U^+(t_{k+1});    \quad {p}_{k+1} \approx \di(]t_k,t_{k+1}]),  \quad {P}_{k+1} \approx \dI(]t_k,t_{k+1}]).
 \end{eqnarray}
Note that the unknown variable $P_{k+1}$ is equivalent to an impulse. This characteristic feature of the Moreau--Jean scheme renders the numerical integration consistent when an impact occurs. Indeed, when the time--step vanishes, a choice of a variable equivalent to a contact force would lead to some unbounded values.

The index set $\mathcal I_1$ results from the time--discretization of the conditional statement in~(\ref{eq:signo-velocity})
\begin{equation}
   \text{if }  g^\alpha(t) \leq 0 \text{ and } U^-(t) <0
\end{equation}
 that allows us to apply the Signorini  condition at the velocity level. In the numerical practice, we choose to define this set by
 \begin{equation}
   \label{eq:index-set1}
   \mathcal I_1 = \{\alpha \in \mathcal I \mid G^\top (q_{k} + h v_{k}) + w \leq 0\text{ and } U_k \leq 0 \}.
 \end{equation}
Other strategies for defining $\mathcal I_1$ can be found in~\citep{Acary2013224}.  The numerical scheme which solves~(\ref{eq:signo-velocity}) enforces in discrete time the Newton law at each time step. Conversely, the constraints $g(t)\geq 0$ are not satisfied. A violation of the constraints, proportional to the time--step, can occur at the activation of the contact, that is, when the contact is closing. The violation may be corrected by a projection technique onto the constraints extending the Gear--Gupta--Leimkuhler approach for DAEs as it has recently been done in~\citep{Acary2013224}.

\subsection{Nonsmooth Newmark and generalized-$\alpha$ scheme}

In \citep{Chen.ea2012,Chen.ea_IJNME2013}, a new family of time--integration schemes has been developed that keep the advantages of the Moreau--Jean in terms of robustness and efficiency  while adding the some key properties of the Newmark~\citep{Newmark59}, the Hilber--Hughes--Taylor (HHT)~\citep{Hughes87} and the generalized-$\alpha$~\citep{Chung.Hulbert93} schemes. The most well--known property of the latter schemes is the  controlled damping of the high--frequency contents of the dynamics. In the linear time invariant dynamics with  unilateral constraints, the new scheme can be written as follows
\begin{subnumcases}{\label{eq:NonsmoothGenerelizedAlphaTS}}
M \dot{\tilde v}_{k+1}  + K q_{k+1} + C {v_{k+1}} = F_{k+1}, \\[1mm]
M w_{k+1} =  G P_{k+1}, \\[1mm]
U_{k+1} = G^\top v_{k+1}, \\[1mm]
\begin{array}{lcl}
  0 \leq U^\alpha_{k+1} + e  U^\alpha_{k} \perp P^\alpha_{k+1}  \geq 0,& \quad&\alpha \in \mathcal I_1, \\[1mm]
  P^\alpha_{k+1}  =0,&\quad& \alpha \in \mathcal I \setminus \mathcal I_1,
\end{array}\label{eq:NonsmoothGenerelizedAlphaTSd} \\[1mm]
\text{with} \notag \\[1mm] 
(1-\alpha_m) a_{k+1} + \alpha_m a_k = (1-\alpha_f) \dot{\tilde v}_{k+1} + \alpha_f  \dot{\tilde v}_{k},\label{eq:NonsmoothGenerelizedAlphaTSe} \\
\tilde q_{k+1} = q_ k + h v_k + h^2 \left(1/2 - \beta\right) a_k + h^2 \beta  a_{k+1}, \\[1mm]
\tilde v_{k+1} = v_k + h (1-\gamma) a_k + h \gamma a_{k+1}, \\[1mm]
v_{k+1} = \tilde v_{k+1} + w_{k+1}, \\[1mm]
q_{k+1} = \tilde q_{k+1} + \frac{h}{2} \,  w_{k+1}.
\end{subnumcases} 

The numerical scheme~(\ref{eq:NonsmoothGenerelizedAlphaTS}) has been designed such that it deals, as the Moreau--Jean scheme with the contact forces and impact through their associated impulses $P_{k+1}$ in a fully implicit way. In this manner, we ensure that the scheme will be consistent when the time--step vanishes if an impact occurs. Furthermore, it also includes a treatment of the unilateral constraint together with the Newton--impact law at the velocity level. This aspect is crucial for the conservation and dissipation properties as we will see in Section~\ref{Sec:Dissipation-Alpha-Scheme}. Finally, the last important property is the second order approximation of the smooth terms given by the generalized--$\alpha$ schemes that allows us to take advantage of the controlled damping of the high-frequency dynamics. It is indeed well--known that one of the main advantages of the $\alpha$--scheme with respect to the Newmark scheme is  the possibility to introduce some damping of the high--frequency dynamics without altering the order.  In the smooth case when the contact is not taken into account, the condition of second order accuracy reads as
\begin{equation}\label{eq:NonsmoothGenerelizedAlphaTS-ordercondition}
 \gamma = \Frac 1 2 + \alpha_f -\alpha_m.
\end{equation}
 The optimal parameters are usually chosen according to the spectral radius at infinity $\rho_\infty \in [0,1]$ such that
\begin{equation}\label{eq:NonsmoothGenerelizedAlphaTS-parameters}
\alpha_m = \Frac{2\rho_\infty-1}{\rho_\infty+1},\quad \alpha_f= \Frac{\rho_\infty}{\rho_\infty+1},\quad \beta = \Frac 1 4 (\gamma+\Frac 1 2 )^2.
\end{equation}
The nonsmooth Newmark algorithm is obtained with $\alpha_m=\alpha_f=0$ and the nonsmooth HHT scheme in the form published in~\citep{Hughes87} is obtained with $\alpha_m=0$ and $\alpha_f \in [0,1/3]$.

\paragraph{Variant of the Moreau--Jean scheme for $2\beta=\gamma =\theta$.} Let us remark that the nonsmooth Newmark algorithm can be reformulated as
\begin{subnumcases}{\label{eq:NonsmoothNewmarkTS}}
M(v_{k+1}-v_{k}) = h (F_{k+\gamma} - K q_{k+\gamma} -C v_{k+\gamma}) + G P_{k+1}\\[1mm]
q_{k+1} = q_ k + h v_k + \Frac 1 2 h M^{-1}[h (F_{k+2\beta} - K q_{k+2\beta} -C v_{k+2\beta}   ) + G P_{k+1})]
\end{subnumcases} 
with~(\ref{eq:NonsmoothGenerelizedAlphaTSd}).  With straightforward manipulations, one obtains if  $2\beta=\gamma =\theta$
\begin{subnumcases}
  {\label{eq:MoreauTS-variant}}
  M(v_{k+1}-v_k) + h   K q_{k+\theta} + h   C v_{k+\theta} - h F_{k+\theta} = p_{k+1} = G P_{k+1},\quad\,\\[1mm] 
  q_{k+1} = q_{k} + h v_{k+1/2}, \quad \\[1mm]
  U_{k+1} = G^\top \, v_{k+1}, \\[1mm]
  \begin{array}{lcl}
    0 \leq U^\alpha_{k+1} + e  U^\alpha_{k} \perp P^\alpha_{k+1}  \geq 0,& \quad&\alpha \in \mathcal I_1, \\[1mm]
    P^\alpha_{k+1}  =0,&\quad& \alpha \in \mathcal I \setminus \mathcal I_1.
  \end{array}
\end{subnumcases} 
This scheme appears as a variant of the Moreau--Jean scheme. We will see in Section~\ref{Sec:NonSmoothNewmark} that these scheme has better dissipation properties than the original one. Note that the original Moreau--Jean scheme cannot be viewed as a special case of the nonsmooth Newmark scheme when $\theta\neq 1/2$.

\subsection{Nonsmooth Krenk--H{\o}gsberg (KH) generalized--$\alpha$ scheme}

In \citep{Krenk.Hogsberg2005} and \citep{Krenk2006}, an alternative collocation method is proposed for the generalized--$\alpha$ scheme
\begin{equation}
  \label{eq:NewAlpha-1}
  \begin{array}{l}
    (1-\alpha_m) [ M a_{k+1}  + C v_{k+1} -  F_{k+1} ]  + \alpha_m [M a_k  + C v_{k} + F_k  ] = (1-\alpha_f) [- K q_{k+1} ]
    + \alpha_
f [ - K q_{k}].
  \end{array}
\end{equation}
This scheme, that will be termed in the sequel as the Krenk--H{\o}gsberg (KH)  generalized--$\alpha$ scheme, differs only from the original generalized--$\alpha$ scheme by the fact that the damping terms and the load terms have the same weight as the inertial term.  In the original generalized--$\alpha$ scheme as it is given in~(\ref{eq:link-Krenk-2}), the weighting of the damping and the load terms follow the stiffness term. The nonsmooth KH generalized--$\alpha$ scheme is obtained by replacing ~(\ref{eq:NonsmoothGenerelizedAlphaTSd}) by (\ref{eq:NewAlpha-1}) in~(\ref{eq:NonsmoothGenerelizedAlphaTS}). Although the order of the method decreases, we will see in the sequel that this scheme has better dissipation properties than the standard generalized--$\alpha$ scheme.
 
The analysis of the KH scheme in~\citep{Krenk.Hogsberg2005} shows that the scheme introduces a slight increase of the frequency response of the mechanical system. In other terms,  the original generalized--$\alpha$ scheme has a slightly improved frequency response. This is mainly related to the alternative choice of weighting of the structural damping term. On the contrary, the choice of weighting of the load term in the KH scheme is superior.

Most of the original $\alpha$--schemes are contained in the KH  generalized--$\alpha$ scheme. Note that in the original paper on the HHT scheme~\citep{Hilber.Hughes.ea77}, the weighting of the load term follows the inertial term. This is the reason why the original HHT scheme in~\citep{Hilber.Hughes.ea77} can be obtained from (\ref{eq:NewAlpha-1}) with $\alpha_m=0$ and $\alpha_f=\alpha$. In the same manner, the $\alpha$--method of \citep{Wood.Bossak.ea81} can be obtained from~(\ref{eq:KrenkHogsberg3}) by choosing $\alpha_f= 0$.


\section{Energy conservation and dissipation properties of Moreau--Jean scheme}
\label{Sec:EnergyAnalysisMoreau}

In this section, we  give a first result on the energy conservation and dissipation of the Moreau--Jean scheme. This result gives a criteria on the parameter $\theta$ which depends on the coefficient of restitution. Let us define the discrete approximation of the work done by the external forces within the step by
  \begin{equation}
    \label{eq:variation-energy70}
    \newW^{\sf ext}_{k+1} = h v_{k+\theta}^\top F_{k+\theta} \approx \displaystyle\int_{t_k}^{t_{k+1}} F v\,\dt,
  \end{equation}
  and the  discrete approximation of the work done by the damping term by
  \begin{equation}
    \label{eq:variation-energy70-damping}
    \newW^{\sf damping}_{k+1} = - h v_{k+\theta}^\top C v_{k+\theta} \approx - \displaystyle\int_{t_k}^{t_{k+1}} v^T C v\,\dt.
  \end{equation}
  We have the following  estimate for the variation of the total mechanical energy.
\begin{lemma}
  \label{Lemma:MoreauTS}
The discrete--time dissipation equality of the Moreau--Jean scheme~(\ref{eq:MoreauTS}) over a time--step $[t_k,t_{k+1}]$ is given by
  \begin{equation}
    \label{eq:variation-energy4}
    \begin{array}{lcl}
      \Delta \mathcal E - \newW^{\sf ext}_{k+1} - \newW^{\sf damping}_{k+1}
      &=& (\Frac 1 2 - \theta)\left[ \|v_{k+1} - v_{k}\|^2_M + \|(q_{k+1} - q_{k})\|^2_K \right]   + U_{k+\theta}^\top P_{k+1}.
    \end{array}
  \end{equation}
\end{lemma}

\paragraph{Proof}

Let us first compute the variation of the energy $\mathcal E$ over the time--step
\begin{equation}
  \label{eq:variation-energy1}
  \begin{array}{lcl}
  \Delta \mathcal E &=& \mathcal E(q_{k+1},v_{k+1}) - \mathcal E(q_{k},v_{k}) \\
  &=& \Frac 1 2 \left[  (v_{k+1}+v_k)^\top M (v_{k+1}-v_k) + (q_{k+1}+q_k)^\top K (q_{k+1}-q_k)    \right],
\end{array}
\end{equation}
since we assume that $M=M^\top$ and $K=K^\top$. Let us also remark that for the $\theta$-method, we have
  \begin{equation}
    \label{eq:variation-energy2}
    \begin{array}{lcl}
      \Frac  1 2 (v_{k+1} + v_{k}) = \Frac 1 h (q_{k+1}-q_{k}) + (\Frac 1 2 - \theta) (v_{k+1} - v_{k}),
    \end{array}
  \end{equation}
and
  \begin{equation}
    \label{eq:variation-energy2bis}
    \begin{array}{lcl}
      v_{k+\theta} = \Frac  1 2 (v_{k+1} + v_{k}) -  (\Frac 1 2 - \theta) (v_{k+1} - v_{k}) .
    \end{array}
  \end{equation}
Using~(\ref{eq:variation-energy2}) and then~(\ref{eq:MoreauTS-motion}), the energy balance~(\ref{eq:variation-energy1}) becomes
\begin{equation}
  \label{eq:variation-energy3}
  \begin{array}{lcl}
   \Delta \mathcal E
     &=&   (\Frac 1 2 - \theta) (v_{k+1} - v_{k})^\top M (v_{k+1}-v_k) + \Frac 1 h (q_{k+1}-q_{k})  M (v_{k+1}-v_k)\\
     & & \quad + \Frac 1 2 \left[  (q_{k+1}+q_k)^\top K (q_{k+1}-q_k)    \right]  \\[3mm]     
     &=&   (\Frac 1 2 - \theta) (v_{k+1} - v_{k})^\top M (v_{k+1}-v_k) + \Frac 1 h (q_{k+1}-q_{k}) \left[- h   K q_{k+\theta} - h Cv_{k+\theta} + h F_{k+\theta} + G P_{k+1}\right]  \\
     &&  +\Frac 1 2 \left[  (q_{k+1}+q_k)^\top K (q_{k+1}-q_k)    \right]  \\[3mm]
     &=&   (\Frac 1 2 - \theta) (v_{k+1} - v_{k})^\top M (v_{k+1}-v_k) +  (\Frac 1 2 - \theta) (q_{k+1} - q_{k})^\top K (q_{k+1}-q_k) \\
     & & -  h v_{k+\theta}^\top C v_{k+\theta}  + h v_{k+\theta}^\top F_{k+\theta} +v_{k+\theta}^\top G P_{k+1}.
  \end{array}
\end{equation}
Using the expression of the norms $\|\cdot\|_M$ and the semi--norm $\|\cdot\|_K$, the expression~(\ref{eq:variation-energy3})  can be easily simplified in
\begin{equation}
  \label{eq:variation-energy4a}
  \begin{array}{lcl}
    \Delta \mathcal E
    &=& (\Frac 1 2 - \theta)\left[ \|v_{k+1} - v_{k}\|^2_M + \|(q_{k+1} - q_{k})\|^2_K \right]  +  h v_{k+\theta}^\top F_{k+\theta} -  h v_{k+\theta}^\top C v_{k+\theta} +v_{k+\theta}^\top G P_{k+1}.
  \end{array}
\end{equation}
Using the definition of the discrete works in (\ref{eq:variation-energy70}) and (\ref{eq:variation-energy70-damping}) and  the fact that $v_{k+\theta}^\top G P_{k+1} = U_{k+\theta}^\top P_{k+1}$ yields the result. \qedtext

\begin{remark} The previous result may be specified at the order $h$.
  Using the following  approximation for function of bounded variations~\citep{acary:inria-00476398}
  \begin{equation}
    \label{eq:variation-energy4bis}
    h v_{k+\theta}^\top F_{k+\theta} - \int_{t_k}^{t_{k+1}}  F(t)v(t) \, \mathrm{d}t=  \mathcal O(h) ,
  \end{equation}
  and 
  \begin{equation}
    \label{eq:variation-energy4ter}
    h v_{k+\theta}^\top C v_{k+\theta} - \int_{t_k}^{t_{k+1}}  v^\top(t) C v(t) \, \mathrm{d}t=  \mathcal O(h) ,
  \end{equation}
  we get for $h$ small enough
  \begin{equation}
    \label{eq:variation-energy40}
    \begin{array}{cl}
      & \displaystyle \Delta \mathcal E  - \int_{t_k}^{t_{k+1}}  Fv \, \mathrm{d}t + \int_{t_k}^{t_{k+1}}  C v \, \mathrm{d}t  = \\
      & \hphantom{ \displaystyle \Delta \mathcal E  - \int_{t_k}^{t_{k+1}}  Fv \, \mathrm{d}t } (\Frac 1 2 - \theta)\left[ \|v_{k+1} - v_{k}\|^2_M + \|(q_{k+1} - q_{k})\|^2_K \right]   + U_{k+\theta}^\top P_{k+1} +  \mathcal O(h).
    \end{array}
  \end{equation}
\end{remark}

Let us give now a first result concerning the dissipation of the Moreau--Jean scheme
\begin{proposition} 
  \label{Prop:Moreau-Jean} The Moreau--Jean scheme  dissipates energy in the sense that 
  \begin{equation}
     \mathcal E(t_{k+1}) - \mathcal E(t_{k})   \leq \newW^{ \sf ext}_{k+1}  + \newW^{\sf damping}_{k+1},
  \end{equation}
  if 
  \begin{equation}
    \label{eq:variation-energy10}
    \Frac 1 2 \leq \theta \leq \Frac 1 {1+e} \leq 1.
  \end{equation}
where $e = \max{e^\alpha, \alpha \in \mathcal I}$. In particular, for $e =0$, we get $\Frac 1 2 \leq\theta \leq 1$ and for $e=1$, we get  $\theta = \Frac 1 2$.
\end{proposition}
In other words,  providing that (\ref{eq:variation-energy10}) is satisfied, the variation of the total mechanical energy of the system is always less than the energy supplied by the external and damping forces.

\paragraph{Proof} Obviously, we have
\begin{equation}
  \label{eq:variation-energy5}
    (\Frac 1 2 - \theta) \left[ \|v_{k+1} - v_{k}\|^2_M + \|(q_{k+1} - q_{k})\|^2_K \right]  \leq 0, \text{ if and only if } \theta \geq 1/2.
\end{equation}
I remains to prove that $U_{k+\theta}^\top  P_{k+1} \leq 0 $. Let us define the following index set of contacts
\begin{equation}
  \label{eq:index-set-impact}
  \begin{array}{lcl}
    \mathcal I^0_1 &=& \{\alpha \in \mathcal I_1 \mid P^\alpha_{k+1} =0\}, \\
  \end{array}
\end{equation}
and its complement in $\mathcal I_1$ is denoted by $\overline{ \mathcal I^0_1 } =\mathcal I_1 \setminus   \mathcal I^0_1 $.
We can therefore develop  $U_{k+\theta}^\top  P_{k+1}  $ as
\begin{equation}
  \label{eq:variation-energy6-1}
  \begin{array}{lcl}
   U_{k+\theta} ^\top P_{k+1} = \displaystyle \sum_{\alpha \in \mathcal I}  U^\alpha_{k+\theta} P^\alpha_{k+1} 
     &=& \displaystyle \sum_{\alpha \in \mathcal I_1}  U^\alpha_{k+\theta} P^\alpha_{k+1} \text{ since  for } \alpha \not\in \mathcal I_1, \quad P^\alpha_{k+1}=0 \\
     &=& \displaystyle \sum_{\alpha \in \overline{ \mathcal I^0_1 }}  (1-\theta(1+e)) U^\alpha_k P^\alpha_{k+1}\\
 \end{array}
\end{equation}
 since  for $\alpha \in \mathcal I^0_1, P^\alpha_{k+1}=0$ and   $\alpha \in \overline{\mathcal I^0_1},  U^\alpha_{k+1}=-e U^\alpha_k$. Since we have  $P_{k+1}^\alpha \geq 0$ and   $U^\alpha_{k} \leq 0$ for all $\alpha \in \mathcal I_1$,  the constraint on $\theta$ is therefore
\begin{equation}
  \label{eq:variation-energy7}
  \theta \leq \Frac 1 {1+e^\alpha} \leq 1, \text{ for all } \alpha \in \mathcal I.
\end{equation}
By combining the constraints on $\theta$ in (\ref{eq:variation-energy5}) and (\ref{eq:variation-energy7}), the result is proved. \qedtext

The following comments can be made on Proposition \ref{Prop:Moreau-Jean}:
  \begin{enumerate}
  \item The variation of energy (\ref{eq:variation-energy4}) may be
    formulated in another form as
    \begin{equation}
      \label{eq:variation-energy11}
      \begin{array}{lcl}
        \Delta \mathcal E - \bar W^{\sf ext}_{k+1}- \newW^{\sf damping}_{k+1}
        &=& (\Frac 1 2 - \theta) \left [\|v_{k+1} - v_{k})\|^2_M +  \|q_{k+1} - q_{k}\|^2_K\right]\\[1mm]
        & & + \Frac 1 2  (v_{k+1}+v_k)^\top G P_{k+1}  - (\Frac 1 2-\theta)  (v_{k+1}-v_k)^\top G P_{k+1} \\[3mm]
        &=& (\Frac 1 2 - \theta) \left [\|v_{k+1} - v_{k})\|^2_M +  \|q_{k+1} - q_{k}\|^2_K  - P^\top_{k+1} (U_{k+1}-U_k)\right]\\[1mm]
        & & +  P^\top_{k+1} U_{k+1/2}.
      \end{array}
    \end{equation}
    The term $ P^\top_{k+1} U_{k+1/2}$ appears as the  discrete dissipated energy at impact. We have also
    \begin{equation}
      \label{eq:variation-energy13}
      \begin{array}{l}
        - P_{k+1}^\top (U_{k+1}-U_k)  =  \sum_{\alpha\in\mathcal I_1}  (1+e) P^\alpha_{k+1} U^\alpha_k \leq 0.
      \end{array}
    \end{equation}
    This alternative form~(\ref{eq:variation-energy11}) shows that the scheme is always dissipative for $\theta = \frac 1 2$.

    \item The bounds (\ref{eq:variation-energy10}) are not sharp since a part
    of the energy potentially generated at impact  $(1-\theta(1+e)) U^\top_k
    P_{k+1}$ is dissipated by the inertial term      $(\Frac 1 2 - \theta) \left [\|v_{k+1} - v_{k})\|^2_M \right]$.

  \item For    $e=0$, the scheme is dissipative for the whole range $\theta \in [\frac 1 2 , 1]$. We can also observe that the trend is somehow opposed the standard property of the $\theta$--method. If the system is more dissipative at contact when $e^\alpha =1$, we have to use the most conservative case with $\theta=1/2$.   
  \end{enumerate}

\section{Dissipation properties of nonsmooth generalized-$\alpha$ scheme}
\label{Sec:Dissipation-Alpha-Scheme}
In this section, the behavior of the nonsmooth generalized-$\alpha$ scheme concerning the energy-conserving or dissipating properties is studied. The proposed method of study  is an extension of the pioneering works of  Krenk and H{\o}gsberg~\citep{Krenk.Hogsberg2005,Krenk2006} on the $\alpha$--schemes. 


\subsection{Principle  of the Krenk--H{\o}gsberg  method and its extension}
\label{Sec:Krenk-Principle}
One of the fundamental properties of the generalized-$\alpha$ scheme is the introduction of a controllable damping of the high--frequency dynamics without altering the second--order accuracy. In~\citep{Krenk.Hogsberg2005,Krenk2006}, this property is studied by explicitly exhibiting and adding a first--order filter and an associated additional state variable in the time--continuous dynamics. Once the augmented time--continuous dynamics is defined, it is shown that the generalized--$\alpha$ scheme applied to the original dynamics is equivalent to the application of the standard Newmark scheme to the augmented dynamics  over two consecutive time--steps and performing a weighting procedure. Hence, the study that have been done for the Newmark scheme in the previous section can be adapted to the generalized-$\alpha$ scheme.

The original Krenk--H{\o}gsberg  method in~\citep{Krenk.Hogsberg2005,Krenk2006} is based on the introduction one additional filter and one additional variable to the original dynamics. The augmented dynamics
\begin{equation}
  \label{eq:Krenk1-ori}
  M a(t) + C v(t) + K q(t) = F(t) + \Frac {\eta}{\nu} [ K z(t) ],
\end{equation}
is introduced, together with the following auxiliary dynamics that filter the previous one
\begin{equation}
  \label{eq:Krenk2-ori}
  \nu h \dot z(t) + z(t) = \nu h \dot q(t), 
\end{equation}
where the time scale of the filter is given by $\nu h$. The parameter $\eta$ is a non--dimensional parameter that permits to tune the effect of the filter on the original dynamics. The dynamics of the first order filter in~(\ref{eq:Krenk2-ori}) is discretized by means of a mid--point rule:
\begin{equation}
  \label{eq:Krenk4-ori}
  \nu (z_{k+1} - z _k) +\frac 1 2 (z_{k+1}+z_{k}) = \nu (q_{k+1}- q_k), 
\end{equation}
Rearranging the terms, we equivalently write (\ref{eq:Krenk4-ori}) as
\begin{equation}
  \label{eq:Krenk5-ori}
  (\frac 1 2 +\nu) z_{k+1}  + (\frac 1 2 -\nu)  z _k = \nu (q_{k+1}- q_k).
\end{equation}
Let us consider now a linear combination of the augmented equation of motion (\ref{eq:Krenk1-ori}) with the weight $(1 /2 + \nu)$ at time $t_{k+1}$ and the weight $(1/ 2 - \nu)$ at time $t_k$:
\begin{equation}
  \label{eq:KrenkHogsberg3}
  \begin{array}{lcl}
    (\Frac 1 2 +\nu) [M a_{k+1} + C v_{k+1} - F_{k+1}] &+& (\Frac 1 2 -\nu) [M a_k + C v_k -F_k]\\
    &=& (\Frac 1 2 +\nu -\eta) [ - K q_{k+1} ]+ (\Frac 1 2 -\nu +\eta) [ - K q_{k} ].
\end{array}
\end{equation}
By choosing the values of $\nu$ and $\eta$ such that
 \begin{equation}
  \label{eq:link-Krenk}
  \begin{array}{l}
  \nu = \Frac 1 2 -\alpha_m \\
  \eta  =  \nu  - \Frac 1 2 + \alpha_f = \alpha_f-\alpha_m,
\end{array}
\end{equation}
 we obtain the KH generalized--$\alpha$ scheme as in~(\ref{eq:NewAlpha-1}).
The standard energetic analysis of the Newmark scheme can be then extended to the KH generalized-$\alpha$ scheme by  adding the following damping force in the energetic  analysis of the Newmark scheme
\begin{equation}
  \label{eq:Krenk10-ori}
  f^A =  \Frac {\eta}{\nu} [K z]
\end{equation}
as it has been shown in~\citep{Krenk2006}. This result will be extended to the nonsmooth case in Section~\ref{Sec:Krenk-Variant}.

 Unfortunately, the previous approach does not longer directly apply to the study of the standard generalized--$\alpha$ scheme since the weighting of the damping and the load terms follow the stiffness term in the method presented by \cite{Chung.Hulbert93}. In the following, we use three additional filters and three associated variables $x,y,z$ to the original dynamics. Let us introduce the augmented dynamics
\begin{equation}
  \label{eq:Krenk1}
  M a(t) + C v(t) + K q(t) = F(t) + \Frac {\eta}{\nu} [ K z(t) +  C x(t) - y(t) ],
\end{equation}
and the following auxiliary dynamics that filter the previous one
\begin{equation}
  \label{eq:Krenk2}
  \begin{array}{lcl}
  \nu h\, \dot z(t) + z(t) &=& \nu h\, \dot q(t), \\
  \nu h\, \dot x(t) + x(t) &=& \nu h\, \dot v(t), \\
  \nu h\, \dot y(t) + y(t) &=& \nu h\, \dot F(t).
\end{array}
\end{equation}
As previously,  let us consider now a linear combination of the augmented equation of motion (\ref{eq:Krenk1}) with the weight $( 1/ 2 + \nu)$ at time $t_{k+1}$
 and the weight $(1/ 2 - \nu)$ at time $t_k$. By choosing the values of $\nu$ and $\eta$ as in~(~\ref{eq:link-Krenk}), the following discretization is obtained
\begin{equation}
  \label{eq:link-Krenk-2}
  \begin{array}{l}
  (1-\alpha_m) M a_{k+1} + \alpha_m M a_k = (1-\alpha_f) [- K q_{k+1} - C v_{k+1} +  F_{k+1}]
+ \alpha_f [ - K q_{k} - C v_{k}+  F_{k}]   .
\end{array}
\end{equation}
The relation (\ref{eq:link-Krenk-2}) is the characteristic relation that defines the generalized-$\alpha$ scheme. Contrary to the work in~\citep{Krenk2006}, the physical meaning of the filters are more difficult to justify, but it enables to retrieve the second order accurate generalized--$\alpha$ scheme developed in~\citep{Chung.Hulbert93}. 

\begin{remark}
  In the sequel, we will assume the ratio $\eta/\nu$ is finite which is not necessarily the case if $\nu=0$. For instance, if $\rho_\infty=1$, we get from~(\ref{eq:NonsmoothGenerelizedAlphaTS-parameters})
  \begin{equation}\label{eq:NonsmoothGenerelizedAlphaTS-specific1}
    \alpha_m = \Frac 1 2,\quad \alpha_f= \Frac 1 2,
  \end{equation}
  that yields
  \begin{equation}\label{eq:NonsmoothGenerelizedAlphaTS-ordercondition-specific1}
    \nu = 0 \text{ and } \eta =0.
  \end{equation}
  We known that the case  $\rho_\infty=1$ corresponds to the case with the minimal damping. The filters whose time-scale vanishes do not act as a filter since we get from (\ref{eq:Krenk2}) that $z(t) = 0$. However, from~(\ref{eq:NonsmoothGenerelizedAlphaTS-parameters}), we obtain also
  \begin{equation}\label{eq:specific10}
    \frac \eta \nu  = \Frac{1/2-\alpha_m}{\alpha_f-\alpha_m} = \Frac{\Frac 1 2-\Frac {2\rho_\infty-1}{\rho_\infty+1}}{\Frac{\rho_\infty}{\rho_\infty+1}  -\Frac {2\rho_\infty-1}{\rho_\infty+1} } = \Frac 2 3,
  \end{equation}
so that the subsequent analysis remains valid.
\end{remark}


 The standard energetic analysis of the Newmark scheme can be then extended to the generalized-$\alpha$ scheme by  adding the following damping force in the energetic  analysis of the Newmark scheme
\begin{equation}
  \label{eq:Krenk10}
  f^A =  \Frac {\eta}{\nu} [K z +  C x -  y].
\end{equation}
Let us define a discrete ``algorithmic energy'' of the form
\begin{equation}
  \label{eq:Krenk15}
  \mathcal H(q,v,a,z) = \mathcal E(q,v) + \Frac {h^2} 4  ( 2 \beta-\gamma) a^\top M  a
  + \Frac {\eta}{2\nu^2} ( \nu - (\gamma - \Frac 1 2)) z ^\top  K z.
\end{equation}
Let us define the discrete approximation of the work done by the external forces within the step by
\begin{equation}
  \label{eq:Alpha-variation-energy70}
  \newW^{\sf ext}_{k+1} = (q_{k+1}-q_{k})^\top  F_{k,\gamma}   \approx \displaystyle\int_{t_k}^{t_{k+1}} F v\,\dt
\end{equation}
and the  discrete approximation of the work done by the damping term by
\begin{equation}
  \label{eq:Alpha-variation-energy70-damping}
  \newW^{\sf damping}_{k+1} = - (q_{k+1}-q_{k})^\top   C v_{k,\gamma}    \approx - \displaystyle\int_{t_k}^{t_{k+1}} v^T C v\,\dt.
\end{equation}

The following result can be obtained.
\begin{lemma}\label{lemma:Krenk-variation-algoenergy}
 The variation of the ``algorithmic'' energy $\Delta \mathcal H$ over a time--step performed by the nonsmooth generalized-$\alpha$ scheme~(\ref{eq:NonsmoothGenerelizedAlphaTS}) is 
  \begin{equation}
    \label{eq:Krenk-variation-algoenergy71}
    \begin{array}{lcl}
      \Delta \mathcal H - \newW^{\sf ext}_{k+1} - \newW^{\sf damping}_{k+1} &+&  (q_{k+1}-q_{k})^\top\Frac{ \eta} {\nu}[y_{k+\gamma} -C x_{k+\gamma}] \\
      &=&   U_{k+1/2}^\top P_{k+1}  + \Frac 1 2 h^2 (\Frac 1 2 - \gamma)  ( 2 \beta-\gamma) \|(a_{k+1}-a_{k})\|^2_M\\
      & & + (\eta + \Frac 1 2 - \gamma) \|q_{k+1}-q_k\|^2_K   + \Frac{\eta}{\nu} ( \gamma  -\nu - \Frac 1 2) \|z_{k+1}-z_k\|^2_K .
    \end{array}
  \end{equation}
\end{lemma}
The proof of this Lemma is given in Appendix~\ref{Ann:Proof1}.

Let us remark that the analysis of the dissipation properties of the nonsmooth generalized--$\alpha$ scheme only differs from the non--impulsive case by the term $U_{k+1/2} ^\top P_{k+1} $. This is mainly the result of the design of the  nonsmooth generalized--$\alpha$ scheme which deals with the nonsmooth terms with a low--order approximation scheme. By the way, a direct analysis of the nonsmooth Newmark scheme can be carried out by directly extending the work of~\cite{Hughes87}. For the sake of space,  the nonsmooth Newmark scheme will be treated as a special case in Section~\ref{Sec:NonSmoothNewmark}. Let us give a result of the sign of $U_{k+1/2} ^\top P_{k+1}$ that appears as the additional term due to the nonsmooth terms in the dynamics.
\begin{lemma}
  \label{lemma:discretework}
  Let us consider that the local velocities $U_{k+1}$ and impulses $P_{k+1}$ satisfies~(\ref{eq:NonsmoothGenerelizedAlphaTSd}). Then the discrete work of the contact forces is negative
  \begin{equation}
    \label{eq:DiscreteWork}
    U_{k+1/2} ^\top P_{k+1} \leq 0
  \end{equation}
\end{lemma}

\paragraph{Proof}
 By introducing the sets of indices as in (\ref{eq:index-set-impact}), 
we have
 \begin{equation}
  \label{eq:Newmark-variation-energy6-1}
  \begin{array}{lcl}
   U_{k+1/2} ^\top P_{k+1} = \displaystyle \sum_{\alpha \in \mathcal I}  U^\alpha_{k+1/2} P^\alpha_{k+1} 
   &=& \displaystyle \sum_{\alpha \in \mathcal I_1}  U^\alpha_{k+1/2} P^\alpha_{k+1} \text{ since  for } \alpha \not\in \mathcal I_1, \quad P^\alpha_{k+1}=0 \\
   &=& \displaystyle \sum_{\alpha \in \overline{ \mathcal I^0_1 }}  \frac 1 2 (1-e) U^\alpha_k P^\alpha_{k+1}
 \end{array}
\end{equation}
 since  for $\alpha \in \mathcal I^0_1, P^\alpha_{k+1}=0$ and   $\alpha \in \overline{\mathcal I^0_1},  U^\alpha_{k+1}=-e U^\alpha_k$.
We conclude that  $U_{k+1/2} ^\top P_{k+1} \leq 0$ since  $P_{k+1}^\alpha \geq 0$ and   $U^\alpha_{k} \leq 0$ for all $\alpha \in \mathcal I_1$.
\qedtext

Lemma~\ref{lemma:Krenk-variation-algoenergy} and Lemma~\ref{lemma:discretework} do not permit to conclude in the general case to the dissipation of the scheme. This is mainly due to the presence of the terms related to $y_{k+\gamma}$ and $ x_{k+\gamma}$ in the left--hand side of (\ref{eq:Krenk-variation-algoenergy71}). Although these terms are only related to the external forces and the damping terms, it seems difficult to expressed them in terms of the original dynamical system. One can only conclude on special cases when the external forces are constant and the damping matrix vanishes. In the next sections, we prefer to focus our effort on the nonsmooth KH generalized-$\alpha$ scheme in Section~\ref{Sec:Krenk-Variant}, on the HHT scheme in Section~\ref{Sec:HHT-scheme} and on the  nonsmooth  Newmark scheme in Section~\ref{Sec:NonSmoothNewmark}.




\subsection{The nonsmooth KH generalized--$\alpha$ scheme case.}
\label{Sec:Krenk-Variant}

The following proposition is a direct application of Lemma~\ref{lemma:Krenk-variation-algoenergy} and Lemma~\ref{lemma:discretework}
\begin{proposition}\label{lemma:KrenkHogsberg}
 The variation of the ``algorithmic'' energy $\Delta \mathcal H$ over a time--step performed by the nonsmooth KH generalized--$\alpha$ scheme~(\ref{eq:NewAlpha-1}) is 
 \begin{equation}
   \label{eq:KrenkHogsberg4}
   \begin{array}{lcl}
     \Delta \mathcal H - \newW^{\sf ext}_{k+1} - \newW^{\sf damping}_{k+1}
     &=&   U_{k+1/2}^\top P_{k+1}  -\Frac 1 2 h^2 (\gamma-\Frac 1 2)  ( 2 \beta-\gamma) \|(a_{k+1}-a_{k})\|^2_M\\
     & &  - (\gamma-\Frac 1 2 - \eta) \|q_{k+1}-q_k\|^2_K   - \Frac{\eta}{\nu} (\nu- \gamma + \Frac 1 2) \|z_{k+1}-z_k\|^2_K .
   \end{array}
 \end{equation}
Moreover, the nonsmooth KH generalized $\alpha$--scheme is stable in the following sense 
\begin{equation}
  \label{eq:KrenkHogsberg5}
  \Delta \mathcal H - \newW^{\sf ext}_{k+1} - \newW^{\sf damping}_{k+1} \leq 0,
\end{equation}
if
\begin{equation}
  \label{eq:KrenkHogsberg6}
  2\beta \geq \gamma \geq \Frac 1 2 \quad\text{ and }\quad 0 \leq \eta \leq \gamma-\Frac 1 2 \leq \nu.
\end{equation}
In terms of $\alpha_m$ and $\alpha_f$, the condition~(\ref{eq:KrenkHogsberg6}) is equivalent to
\begin{equation}
  \label{eq:KrenkHogsberg7}
  2\beta \geq \gamma \geq \Frac 1 2 \quad\text{ and }\quad 0 \leq \alpha_f-\alpha_m \leq \gamma-\Frac 1 2 \leq \Frac 1 2 - \alpha_m.
\end{equation}
\end{proposition}

\paragraph{Proof:}
The proof of the equation~(\ref{eq:KrenkHogsberg4}) follows exactly the same lines as the proof of Lemma~\ref{lemma:Krenk-variation-algoenergy} by cancelling the term $y_{k+\gamma}$ and $x_{k+\gamma}$. The inequality~(\ref{eq:KrenkHogsberg5}) is directly obtained with the conditions~(\ref{eq:KrenkHogsberg6}) and the fact that $U^\top_{k+1/2}P_{k+1} \leq 0$ comes from Lemma~\ref{lemma:discretework}. 
The equivalent form of the condition in~(\ref{eq:KrenkHogsberg7}) is obtained with the help of (\ref{eq:link-Krenk}).
\qedtext

Note that with the second order accuracy condition~(\ref{eq:NonsmoothGenerelizedAlphaTS-ordercondition}) $\gamma=1/2 + \alpha_f-\alpha_m$, the condition~(\ref{eq:KrenkHogsberg7}) simplifies in
\begin{equation}
  \label{eq:NewAlpha-cond}
  2\beta \geq \gamma \geq \Frac 1 2 \quad\text{ and }\quad 0 \leq \alpha_f-\alpha_m \leq \Frac 1 2-\alpha_m.
\end{equation}
The nonsmooth KH generalized--$\alpha$ scheme appears as an interesting alternative for the computation of the linear elastodynamics of mechanical system with unilateral contact and impact.


\subsection{The nonsmooth HHT case}
\label{Sec:HHT-scheme}

With the special choice $\alpha_m=0$, the generalized-$\alpha$ scheme reduces to the HHT scheme in the form presented in~\citep{Hughes87}. The equivalent filter parameters are $\nu=1/2$ and $\eta = \alpha_f$ for the HHT scheme. For the sake of simplicity, the parameter $\alpha_f$ will be denoted as $\alpha:=\alpha_f$. The HHT scheme is given by
\begin{equation}
  \label{eq:Krenk9}
   M a_{k+1}  + (1-\alpha) [ K q_{k+1} + C v_{k+1}] +\alpha  [K q_{k} + C v_{k}] =  (1 -\alpha)  F_{k+1}+ \alpha  F_{k}.
\end{equation}
The application of Lemma~\ref{lemma:Krenk-variation-algoenergy} in this context leads to the following definition of the approximation of works as follows:
\begin{equation}
  \label{eq:HHT-variation-energy70}
  \newW^{\sf ext}_{k+1} = (q_{k+1}-q_{k})^\top \left[ (1-\alpha) F_{k,\gamma} +\alpha F_{k-1,\gamma} \right] \approx \displaystyle\int_{t_k}^{t_{k+1}} F v\,\dt
\end{equation}
and
\begin{equation}
  \label{eq:HHT-variation-energy70-damping}
  \newW^{\sf damping}_{k+1} = - (q_{k+1}-q_{k})^\top   C \left[ (1-\alpha) v_{k,\gamma} + \alpha v_{k-1,\gamma} \right]    \approx - \displaystyle\int_{t_k}^{t_{k+1}} v^T C v\,\dt.
\end{equation}
The following result is straightforwardly derived as a consequence of Lemma~\ref{lemma:Krenk-variation-algoenergy}.
\begin{proposition}\label{lemma:HHT-Krenk-variation-algoenergy}
 The variation of the ``algorithmic'' energy $\Delta \mathcal H$ over a time--step performed by the nonsmooth HHT scheme (scheme~(\ref{eq:NonsmoothGenerelizedAlphaTS}) with  $\alpha_m=0$ and  $\alpha_f=\alpha$) is 
  \begin{equation}
    \label{eq:Krenk-variation-algoenergy7100}
    \begin{array}{lcl}
      \Delta \mathcal H - \newW^{\sf ext}_{k+1} - \newW^{\sf damping}_{k+1} 
      &=&   U_{k+1/2}^\top P_{k+1}   -\Frac 1 2 h^2 (\gamma-\Frac 1 2)  ( 2 \beta-\gamma) \|(a_{k+1}-a_{k})\|^2_M\\
      &&  - (\gamma -\Frac 1 2 - \alpha) \|q_{k+1}-q_k\|^2_K - 2 \alpha (1 - \gamma) \|z_{k+1}-z_k\|^2_K .
    \end{array}
  \end{equation}
Moreover, the nonsmooth HHT scheme dissipates the  ``algorithmic'' energy $\mathcal H$ in the following sense 
\begin{equation}
  \label{eq:Krenk-variation-algoenergy7101}
  \Delta \mathcal H - \newW^{\sf ext}_{k+1} - \newW^{\sf damping}_{k+1} \leq 0,
\end{equation}
if
\begin{equation}
  \label{eq:Krenk-variation-algoenergy7102}
 2\beta \geq \gamma \geq \Frac 1 2 \quad\text{ and }\quad 0 \leq \alpha \leq \gamma - \Frac 1 2  \leq \Frac 1 2. 
\end{equation}
\end{proposition}

\paragraph{Proof:}
The right--hand side of (\ref{eq:Krenk-variation-algoenergy7100}) is directly obtained from Lemma~\ref{lemma:Krenk-variation-algoenergy} by writing $\nu = 1/2$ and $\eta=\alpha$. We have still to simplify the term related to the external and damping forces. With $\nu = 1 /2$, the discretization of the filters amounts to solving
\begin{equation}
  \label{eq:HHT-proof1}
    \begin{array}{rcl}
     \frac 1 2 (x_{k+1} - x _k) +\frac 1 2(x_{k+1}+x_{k}) &=& \frac 1 2 (v_{k+1}- v_k), \\[2mm]
     \frac 1 2 (y_{k+1} - y _k) +\frac 1 2 (y_{k+1}+y_{k}) &=& \frac 1  2 (F_{k+1}- F_k) .     
    \end{array}
\end{equation}
Simplification yields
\begin{equation}
  \label{eq:HHT-proof2}
   \begin{array}{rcl}
     2 x_{k+1} &=& v_{k+1}-v_k, \\
     2 y_{k+1} &=& F_{k+1}-F_k .
   \end{array}
 \end{equation}
Therefore, we obtained for the additional terms
\begin{equation}
  \label{eq:HHT-proof3}
  \begin{array}{lcl}
  F_{k,\gamma} - \Frac {\eta}{\nu} y_{k,\gamma} &=&   F_{k,\gamma} - 2 \alpha y_{k,\gamma} \\
  &=& \gamma F_{k+1} - (1-\gamma) F_k - \alpha \left[ \gamma (F_{k+1} - F_k) + (1-\gamma) (F_k - F_{k-1})  \right]\\
  &= & (1-\alpha) \left[ \gamma F_{k+1} + (1-\gamma) F_k  \right] +  \alpha \left[ \gamma F_{k} + (1-\gamma) F_{k-1} \right]\\
  &=&  (1-\alpha) F_{k,\gamma} +  \alpha F_{k-1,\gamma}.
\end{array}
\end{equation}
By applying the same manipulations to the term involving $C x_{k+\gamma}$ the result is obtained.

The inequality (\ref{eq:Krenk-variation-algoenergy7101}) is straightforwardly obtained thanks to the positiveness of the quadratic terms when the conditions~(\ref{eq:Krenk-variation-algoenergy7102}) are applied to (\ref{eq:Krenk-variation-algoenergy7100}). Since the remaining term $U_{k+1/2}^\top P_{k+1}$ from Lemma~\ref{lemma:discretework}, the proof is completed. 
\qedtext

Note that with the second order accuracy condition~(\ref{eq:NonsmoothGenerelizedAlphaTS-ordercondition}) $\gamma=1/2 + \alpha$, the condition simplifies in
\begin{equation}
  \label{eq:HHT-cond}
  2\beta \geq \gamma \geq \Frac 1 2 \quad\text{ and }\quad 0 \leq \alpha \leq \Frac 1 2.
\end{equation}

\subsection{The nonsmooth Newmark case}
\label{Sec:NonSmoothNewmark}
With the Newmark scheme ($\alpha_m=\alpha_f=0$), the value of the parameters are $\nu=0, \eta=1/2$. The algorithmic energy reduces to  
\begin{equation}
  \label{eq:Krenk2-Newmark}
  \mathcal H(q,v,a) = \mathcal E(q,v) + \Frac {h^2} 4  ( 2 \beta-\gamma) a^\top M  a.
\end{equation}
Although there is no direct mechanical interpretation of this quantity, it allows one to conclude on the dissipation property of the scheme since $\mathcal K$ is a semi--norm for $2\beta \geq \gamma$. Let us remark that we retrieve the algorithmic energy introduced by~\cite{Hughes1977}.  The following result can be obtained.
\begin{proposition}\label{lemma:Newmark-Krenk-variation-algoenergy}
 The variation of the ``algorithmic'' energy $\Delta \mathcal H$ over a time--step performed by the nonsmooth Newmark scheme~(\ref{eq:NonsmoothNewmarkTS}) is 
  \begin{equation}
    \label{eq:Newmark-Krenk-variation-algoenergy71}
    \begin{array}{lcl}
      \Delta \mathcal H - \newW^{\sf ext}_{k+1} - \newW^{\sf damping}_{k+1}
      &=  &  (\Frac 1 2 -\gamma) \left [ \|q_{k+1}-q_k\|^2_K  +\Frac h 2  ( 2 \beta-\gamma) \|(a_{k+1}-a_{k})\|^2_M\right]+  U_{k+1/2}^\top  P_{k+1} .
    \end{array}
  \end{equation}
  Moreover, the nonsmooth Newmark scheme dissipates the ``algorithmic'' energy $\mathcal H$  in the following sense
  \begin{equation}
    \label{eq:Krenk-variation-algoenergy80}
    \begin{array}{lcl}
      \Delta \mathcal H - \newW^{\sf ext}_{k+1}  - \newW^{\sf damping}_{k+1} \leq 0,
    \end{array}
  \end{equation}
  for 
  \begin{equation}
    \label{eq:Krenk-variation-algoenergy81}
    2 \beta \geq \gamma \geq \Frac 1 2.
  \end{equation}
\end{proposition}

\paragraph{Proof:}
The relation (\ref{eq:Newmark-Krenk-variation-algoenergy71}) is direct application of Lemma~\ref{lemma:Krenk-variation-algoenergy} with $\nu=0, \eta=1/2$.
 The inequality (\ref{eq:Krenk-variation-algoenergy80}) is straightforward to obtain. For  $2 \beta \geq \gamma \geq \Frac 1 2$, the first of the right-hand side of (\ref{eq:Newmark-Krenk-variation-algoenergy71}) is non positive and the fact that $U^\top_{k+1/2}P_{k+1} \leq 0$ comes from Lemma~\ref{lemma:discretework}.
\qedtext\\

The very interesting fact in the nonsmooth Newmark scheme is that the standard dissipation properties of  the original Newmark scheme are conserved when unilateral contacts and impacts are included. It is well-known that considering the unilateral contact at the position level and a semi-implicit rule for the reaction forces leads to blow-up in energy~\citep{Chawla.Laursen98} and an ad-hoc restitution rule at contact (see Figure 4 in~\citep{Chen.ea_IJNME2013}). With the nonsmooth Newmark scheme, the energy balance is equivalent to the unconstrained case. The following result is a corollary of Proposition~\ref{lemma:Newmark-Krenk-variation-algoenergy} that specifies the properties of the nonsmooth Newmark scheme for some particular choices of the parameters $\gamma$ and $\beta$.

\begin{corollary} For the following specific values of the parameters $\gamma$ and $\beta$, the variation of the ``algorithmic'' energy $\Delta \mathcal H$ over a time--step performed by the nonsmooth Newmark scheme~(\ref{eq:NonsmoothNewmarkTS}) can be specified as follows
  \begin{itemize}
  \item For $\gamma = \Frac 1 2 $, we get
    \begin{equation}
      \label{eq:Newmark-Krenk-variation-algoenergy71-spec1}
      \Delta \mathcal H - \newW^{\sf ext}_{k+1} - \newW^{\sf damping}_{k+1} =    U_{k+1/2}^\top  P_{k+1} .
    \end{equation}
    The algorithmic energy $\mathcal K$ only changes  at each time--step by the amount of work that is done by the dissipative effects in the system and the work supplied by the external forces.
  \item For $\gamma=2\beta \geq 1/2$, $\mathcal E(q,v) = \mathcal H(q,v,a)$ and we have
    \begin{equation}
      \label{eq:Newmark-variation-energy-spec2}
      \begin{array}{lcl}
        \Delta \mathcal E -\newW^{\sf ext}_{k+1}  -\newW^{\sf damping}_{k+1}      &=&  (\Frac 1 2-\gamma) \| (q_{k+1}-q_{k } ) \|_K^2 + P^\top_{k+1} U_{k+1/2} 
      \end{array}
    \end{equation}
    Since $ P^\top_{k+1} U_{k+1/2} \leq 0$, the scheme always dissipates the mechanical energy of the system.
  \item  For $\gamma = 2\beta= \Frac 1 2 $, we get
    \begin{equation}
      \label{eq:Newmark-variation-energy-spec3}
      \begin{array}{lcl}
        \Delta \mathcal E -\newW^{\sf ext}_{k+1}  -\newW^{\sf damping}_{k+1}      &=&   P^\top_{k+1} U_{k+1/2} 
      \end{array}
    \end{equation}
     The total mechanical energy $\mathcal E$ only changes at each time--step by the amount of work that is done by the dissipative effects in the system and the work supplied by the external forces.
  \end{itemize}
\end{corollary}

 The results stated in~(\ref{eq:Newmark-variation-energy-spec2}) and~(\ref{eq:Newmark-variation-energy-spec3}) are better than for the standard Moreau--Jean scheme with the $\theta$-method since there is no condition on $\gamma$. The results apply  to the variant of the Moreau-Jean scheme~(\ref{eq:MoreauTS-variant}).


           
\section{Conclusion}
\label{Sec:Conclusion}
The main results of this paper and some perspectives may be summarized now:

\begin{itemize}
\item A study of the energy conservation and dissipation of the Moreau--Jean scheme has been carried out. Under suitable assumptions on the numerical parameter $\theta$ with respect to the coefficients of restitution, we show that the Moreau--Jean scheme dissipates the total mechanical energy. For $\theta=1/2$, the dissipated energy is only due to the discrete works of the damping, external forces and the impact impulses. Moreover, a simple variant of the Moreau--Jean scheme that always dissipates energy has been proposed. This latter method removes the condition on $\theta$ with respect to the coefficients of restitution.
\item As we said in the introduction, there is a lot of difficulties to establish energy conservation or dissipation results in the context of Newmark schemes with unilateral constraints on the generalized coordinates. One of the difficulties raises from the fact that for discrete (or space--discretized) mechanical systems, we have to introduce the notion of coefficient of restitution, or at least to give a law that defines the velocity after the impact. The second difficulty is related to higher order approximation of nonsmooth terms like impulses. With the new schemes developed in~\citep{Chen.ea2012,Chen.ea_IJNME2013}, these two difficulties are overcome. This allows us to extend the known results on the energy conservation and dissipation properties to the nonsmooth case with unilateral contacts and impacts. In particular, we prove that the Newmark scheme is able to conserve or dissipates an algorithmic energy $\mathcal H(q,v,a)$ that depends on the coordinates, the velocities and the acceleration. For HHT scheme, we are also to prove the same kind of results with an extended algorithmic energy $\mathcal H(q,v,a,z)$ based on the introduction of additional filters. For the generalized--$\alpha$ scheme, we obtain equivalent results  as in the smooth case. However, this cannot allow to easily conclude. Nevertheless, the contribution of the contact terms in the discrete energy balance~(\ref{eq:Krenk-variation-algoenergy71}), that is $U^\top_{k+1/2}P_{k+1}$, is identical for the Newmark and the HHT scheme. We can infer that this result combined with the classical result of~\cite{Arnold.Bruls2007,Jay.Negrut2007,Lunk.Simeon2006} on stability of the generalized-$\alpha$ when unilateral constraints remain closed should imply the stability of the scheme in practice. Nevertheless, we extend a variant of the generalized-$\alpha$ scheme due to \cite{Krenk.Hogsberg2005} to the contact case. For this latter case, we obtain the same results as in the smoooth case. 

\item The dissipation properties of the schemes studied in this paper allows us to conclude to the boundedness of the total mechanical energy of the system $\mathcal E(q,v)$ for the Moreau--Jean scheme and to the boundedness the algorithmic energy $\mathcal H(q,v,a,z)$ for the Newmark and the $\alpha$--schemes by adding some standard assumption on the works of the external forces.
Let us remark that the total mechanical energy $\mathcal E(q,v)$ and the algorithmic energy $\mathcal H(q,v,a,z)$ are positive semi--definite functions if we assume that $M$ is positive definite,  $K$ is positive semi-definite and the conditions of dissipation of the schemes are satisfied. The boundedness of theses energies implies the boundedness of the discrete velocities $v_k$ and acceleration $a_k$ which in some sense guarantee the stability of the scheme in the numerical practice. The boundedness of the discrete generalized coordinates can also easily conclude if the stiffness matrix is positive definite, or by simply inspecting the relations that relate the discrete generalized coordinates to the discrete generalized velocities and accelerations. In the smooth case (Lipschitz ordinary differential equation),  boundedness properties are equivalent to the stability of the linear multi--step time integration method thanks to the theory developed by Dalhquist, Lax and Richtmyer~\citep{Dahlquist1985,Lax.Richtmyer1956,Richtmyer.Morton1967}. Classical approaches to study the stability of numerical scheme for ordinary differential equations trough its amplification matrix or the roots of the stability function were also extended with success to the case of differential algebraic equations~\citep{Cardona.Geradin1989,Ascher.Petzold1998,Hairer.Wanner1996}.  Unfortunately, in our case, we cannot directly conclude to the stability of the schemes by any of these methods since we are in a nonsmooth case and the numerical scheme does not reduce to a linear--multistep formula. Nevertheless, it can be a first step to prove the convergence of the scheme.

\item The first perspective that can be drawn for this work is the possible adaptation of the results to the case with Coulomb's friction. One of the main difficulties is that even in the continuous time--case the modeling of Coulomb's friction at the impulse level together with an impact law does not lead necessarily to a dissipative system (see \citep{Leine.vandeWouw2008} for details). The second difficulty should be to prove that the  numerical scheme correctly discretize the friction such that it still dissipates energy. The second perspective is the extension towards to the nonlinear case as it has been done some special case of the mid-point rule in~\citep{Hauret.LeTallec2006,Ayyad.ea2009}.

\end{itemize}

\bibliographystyle{plainnat}
\bibliography{./biblio/String,biblio,./biblio/NonSmooth,./biblio/Math,./biblio/Multibody,./biblio/Fem,./biblio/Dae,./biblio/Meca,./biblio/AnaNum,./biblio/Math-Impact,./biblio/Cp,./biblio/Contact,./biblio/Control}
\appendix

\section{Proof of Lemma~\ref{lemma:Krenk-variation-algoenergy}}
\label{Ann:Proof1}

To prove the result, we start from the dissipation analysis of the Newmark scheme. Following the method introduced in~\citep{Hughes1977} and fully developed in~\citep{Krenk2006}, we define an intermediate discrete ``algorithmic energy'' of the form
\begin{equation}
  \label{eq:Krenk2-bis}
  \mathcal K(q,v,a) = \mathcal E(q,v) + \Frac {h^2} 4  ( 2 \beta-\gamma) a^\top M  a.
\end{equation}
Although there is no direct mechanical interpretation of this quantity, it allows one to conclude on the dissipation property of the scheme since $\mathcal K$ is a semi--norm for $2\beta \geq \gamma$. From the definition  of the nonsmooth Newmark scheme, we get
 \begin{equation}
   \label{eq:Newmark-variation-energy3}
   q_{k+1}-q_{k} =  \Frac h 2 (v_{k+1}+v_{k}) +    \Frac {h^2} 2  ( 2 \beta-\gamma) (a_{k+1}-a_{k})  .
 \end{equation}
and
\begin{equation}
   \label{eq:Hughes3}
   v_{k+1}-v_k =   \Frac h 2 (a_{k+1}+a_k)   + h (\gamma-\Frac 1 2) (a_{k+1}-a_k) +  w_{k+1}.
 \end{equation}
Using (\ref{eq:Newmark-variation-energy3}) and (\ref{eq:Hughes3}), another formulation of $\Delta \mathcal E$ in~(\ref{eq:variation-energy1}) is given by
 \begin{equation}
   \label{eq:Newmark-Krenk0}
   \begin{array}{lcl}
     2 \Delta \mathcal E 
     &=& (v_{k+1}+v_k)^\top M \left[ \Frac h 2 (a_{k+1}+a_k)+ h(\gamma-\Frac 1 2)(a_{k+1}-a_k)+  w_{k+1}\right]\\
     && +  (q_{k+1}+q_{k})^\top K \left[  \Frac h 2 (v_{k+1}+v_{k}) +    \Frac {h^2} 2  ( 2 \beta-\gamma) (a_{k+1}-a_{k}) \right]\\
   \end{array}
 \end{equation}
Using the definition of the nonsmooth Newmark scheme, the relation~(\ref{eq:Newmark-Krenk0}) can be developed as in
 \begin{equation}
   \label{eq:Newmark-Krenk1}
   \begin{array}{lcl}
   2 \Delta \mathcal E 
   &=& (v_{k+1}+v_k)^\top  \left[ \Frac h 2 \big(-K(q_{k+1}+q_k) -C (v_{k+1}+v_k)+ F_{k+1}+F_k \big)   + h M(\gamma-\Frac 1 2)(a_{k+1}-a_k)\right] \\
   &+& 2 U_{k+1/2}^\top P_{k+1} +  (q_{k+1}+q_{k})^\top K \left[  \Frac h 2 (v_{k+1}+v_{k}) +    \Frac {h^2} 2  ( 2 \beta-\gamma) (a_{k+1}-a_{k}) \right]\\
   \\
   &=& (v_{k+1}+v_k) ^\top  \left[ \Frac h 2 ( F_{k+1}+F_k ) - \Frac h 2 C ( v_{k+1}+v_k )  + h M (\gamma-\Frac 1 2)(a_{k+1}-a_k)\right] + 2 U_{k+1/2}^\top P_{k+1}\\
   & & +  \Frac {h^2} 2  ( 2 \beta-\gamma)  (q_{k+1}+q_{k})^\top K(a_{k+1}-a_{k}).
   \\
 \end{array}
\end{equation}
Using again~(\ref{eq:Hughes3}) in~(\ref{eq:Newmark-Krenk1}), we get with some manipulations
\begin{equation}
   \label{eq:Newmark-Krenk2}
   \begin{array}{lcl}
   2 \Delta \mathcal E 
  &=&  ( F_{k+1}+F_k -  C ( v_{k+1}+v_k )   )^\top  \left[    q_{k+1}-q_{k} -   \Frac {h^2} 2  ( 2 \beta-\gamma) (a_{k+1}-a_{k})     \right] + 2 U_{k+1/2}^\top P_{k+1} \\
  & &  + (v_{k+1}+v_k)^\top M \left[h(\gamma-\Frac 1 2)(a_{k+1}-a_k)\right] +   \Frac {h^2} 2  ( 2 \beta-\gamma)  (q_{k+1}+q_{k})^\top K(a_{k+1}-a_{k}) \\
  &=& ( F_{k+1}+F_k -  C ( v_{k+1}+v_k )  )^\top  (q_{k+1}-q_{k}) + 2 U_{k+1/2}^\top P_{k+1}  \\
  & & +   \Frac {h^2} 2  ( 2 \beta-\gamma) (a_{k+1}-a_{k})^\top  \left[ F_{k+1}+F_k - K(q_{k+1}+q_{k}) -  C ( v_{k+1}+v_k )  \right]\\
  & & + h(\gamma-\Frac 1 2)(v_{k+1}+v_k)^\top M (a_{k+1}-a_k) \\
  &=& ( F_{k+1}+F_k - C ( v_{k+1}+v_k )   )^\top  (q_{k+1}-q_{k})  + 2 U_{k+1/2}^\top P_{k+1} \\
  & & +   \Frac {h^2} 2  ( 2 \beta-\gamma) (a_{k+1}-a_{k})^\top  M  (a_{k+1}+a_{k})+  h(\gamma-\Frac 1 2)(v_{k+1}+v_k)^\top M (a_{k+1}-a_k).
 \end{array}
 \end{equation}
Since the term $ (a_{k+1}-a_{k})^\top  M  (a_{k+1}+a_{k}) $ appears in the last line of (\ref{eq:Newmark-Krenk2}), the computation of $ \Delta \mathcal K $ can be developed as follows
 \begin{equation}
   \label{eq:Newmark-Krenk3}
   \begin{array}{lcl}
   2 \Delta \mathcal K 
   &=& ( F_{k+1}+F_k  - C ( v_{k+1}+v_k ) )^\top  (q_{k+1}-q_{k})  +   2 U_{k+1/2}^\top P_{k+1}  \\
   && + h(\gamma-\Frac 1 2)(v_{k+1}+v_k)^\top M (a_{k+1}-a_k) \\
   &=& ( F_{k+1}+F_k  - C ( v_{k+1}+v_k ) )^\top  (q_{k+1}-q_{k})  +   2 U_{k+1/2}^\top P_{k+1} \\
   & & + (\gamma-\Frac 1 2)(a_{k+1}-a_k)^\top M\left[   2(    q_{k+1}-q_{k}) -  h^2  ( 2 \beta-\gamma) (a_{k+1}-a_{k})    \right] \\
   &=& ( F_{k+1}+F_k  - C ( v_{k+1}+v_k )  )^\top  (q_{k+1}-q_{k})  +  2 U_{k+1/2}^\top P_{k+1}    \\
   & &  - (\gamma-\Frac 1 2) \left[  h^2  ( 2 \beta-\gamma)(a_{k+1}-a_k)^\top  M (a_{k+1}-a_{k})    \right] \\
   & &   +2(\gamma-\Frac 1 2)(q_{k+1}-q_k)^\top   \left[     - K(q_{k+1}-q_{k}) - C(v_{k+1}-v_{k})  + F_{k+1}-F_k\right] \\
   &=&  2   (q_{k+1}-q_{k})^\top F_{k+\gamma} - 2 (q_{k+1}-q_{k})^\top C v_{k+\gamma}   +  2 U_{k+1/2}^\top P_{k+1}  \\
   & & -  h^2 (\gamma-\Frac 1 2)  ( 2 \beta-\gamma) \|(a_{k+1}-a_{k})\|^2_M\\
   & & - 2(\gamma-\Frac 1 2) \|q_{k+1}-q_k\|^2_K
  \\
 \end{array}
 \end{equation}
Substituting the definition of  $\newW^{\sf ext}_{k+1}$ and $\newW^{\sf damping}_{k+1}$ in (\ref{eq:Newmark-Krenk3}) yields 
 \begin{equation}
    \label{eq:Newmark-Krenk-variation-algoenergy71-bis}
    \begin{array}{lcl}
      \Delta \mathcal K - \newW^{\sf ext}_{k+1} - \newW^{\sf damping}_{k+1}
      &=  &  (\Frac 1 2 -\gamma) \left [ \|q_{k+1}-q_k\|^2_K  +\Frac h 2  ( 2 \beta-\gamma) \|(a_{k+1}-a_{k})\|^2_M\right]+  U_{k+1/2}^\top  P_{k+1} .
    \end{array}
  \end{equation}
From~(\ref{eq:Alpha-variation-energy70}), the additional damping force (\ref{eq:Krenk10})  generates an additional term in~(\ref{eq:Newmark-Krenk-variation-algoenergy71-bis}) given by
\begin{equation}
  \label{eq:Krenk11}
  \begin{array}{lcl}
    (q_{k+1}-q_k)^\top  f^A_{k+\gamma} 
    &=& \Frac {\eta}{\nu}(q_{k+1}-q_k)^\top [K z_{k+\gamma} +  C x_{k+\gamma} - y_{k+\gamma}  ]
  \end{array}
\end{equation}

For the sake of simplicity, we consider for a moment that $C=0$. The damping term will be discussed later. Developing the additional term
(\ref{eq:Krenk11}) using the rule
\begin{equation}
  z_{k+\gamma} = \frac 1 2 (z_{k+1}+z_k) + (\gamma-\frac 1 2 ) (z_{k+1}-z_k),\label{eq:Krenk3000}
\end{equation}
  we obtain
\begin{equation}
  \label{eq:Krenk12}
  \begin{array}{lcl}
    (q_{k+1}-q_k)^\top  f^A_{k+\gamma} 
    &=& \Frac {\eta}{\nu}(q_{k+1}-q_k)^\top K[ \Frac 1 2 (z_{k+1}+z_k) + (\gamma - \Frac 1 2 )(z_{k+1}-z_k) ] \\ 
    & &- \Frac {\eta}{\nu}(q_{k+1}-q_k)^\top y_{k+\gamma}\\
    &=& \Frac {\eta}{\nu}(q_{k+1}-q_k)^\top K[ \nu (q_{k+1}-q_k) + (\gamma - \Frac 1 2 -\nu)(z_{k+1}-z_k) ]\\
    & & - \Frac {\eta}{\nu}(q_{k+1}-q_k)^\top y_{k+\gamma}\\
    &=&  \eta \|q_{k+1}-q_k\|^2_K + \Frac {\eta}{\nu} (\gamma - \Frac 1 2 -\nu) (q_{k+1}-q_k)^\top K (z_{k+1}-z_k) \\
    & & - \Frac {\eta}{\nu}(q_{k+1}-q_k)^\top y_{k+\gamma} .  \\
  \end{array}
\end{equation}
Let us now focus on the second term in  the right hand side of (\ref{eq:Krenk12}). Using
\begin{equation}
  \label{eq:Krenk3001}
   q_{k+1}-q_k =  z_{k+1}-z_k +\Frac 1 {2 \nu}(z_{k+1}+z_k),
\end{equation}
we obtain
\begin{equation}
  \label{eq:Krenk13}
  \begin{array}{lcl}
    \Frac {\eta}{\nu} (\gamma - \Frac 1 2 -\nu) (q_{k+1}-q_k)^\top K (z_{k+1}-z_k) 
    &=&  \Frac {\eta}{\nu} (\gamma - \Frac 1 2 -\nu) (z_{k+1}-z_k +\Frac 1 {2 \nu}(z_{k+1}+z_k) )^\top       K (z_{k+1}-z_k) \\
    &=&  \Frac {\eta}{\nu} (\gamma - \Frac 1 2 -\nu) \|z_{k+1}-z_k\|^2_K  \\
    & & +\Frac {\eta}{2\nu^2} (\gamma - \Frac 1 2 -\nu) (z_{k+1}+z_k)^\top  K (z_{k+1}-z_k) .
  \end{array}
\end{equation}
Let us restart from (\ref{eq:Newmark-Krenk3}) with the additional term $(q_{k+1}-q_k)^\top  f^A_{k+\gamma}$ , we get
\begin{equation}
  \label{eq:Krenk14}
  \begin{array}{lcl}
    2 \Delta \mathcal K 
    &=&  2 F_{k+\gamma}^\top  (q_{k+1}-q_{k})  - 2 \Frac {\eta}{\nu}(q_{k+1}-q_k)^\top y_{k+\gamma}  +  2 U_{k+1/2}^\top P_{k+1}  \\    
    & &  - h^2 (\gamma-\Frac 1 2)  ( 2 \beta-\gamma) \|(a_{k+1}-a_{k})\|^2_M\\
    & &  - 2(\gamma-\Frac 1 2 - \eta) \|q_{k+1}-q_k\|^2_K  + 2 \Frac {\eta}{\nu} (\gamma - \Frac 1 2 -\nu) \|z_{k+1}-z_k\|^2_K  \\
    & &  + \Frac {\eta}{\nu^2} (\gamma - \Frac 1 2 -\nu) (z_{k+1}+z_k )^\top       K (z_{k+1}-z_k) .
  \end{array}
\end{equation}
With the definition of $\mathcal H$ in~(\ref{eq:Krenk15}), we can simplify (\ref{eq:Krenk14}) in
\begin{equation}
  \label{eq:Krenk16}
  \begin{array}{lcl}
    2 \Delta \mathcal H
    &=& 2 (q_{k+1}-q_{k})^\top [F_{k+\gamma}   -  \Frac {\eta}{\nu} y_{k+\gamma}]  +  2 U_{k+1/2}^\top P_{k+1}  \\
    & & - h^2 (\gamma-\Frac 1 2)  ( 2 \beta-\gamma) \|(a_{k+1}-a_{k})\|^2_M\\
    & & - 2(\gamma-\Frac 1 2 - \eta) \|q_{k+1}-q_k\|^2_K \\
    & & - 2  \Frac {\eta}{\nu} (\nu- \gamma + \Frac 1 2) \|z_{k+1}-z_k\|^2_K  \\
  \end{array}
\end{equation}
The term depending the damping matrix can be included by adding to the external forces a damping force $-Cv$ and the additional term $\frac {\eta}{\nu} x_{k+\gamma}$. Using the definition of the discrete works in~(\ref{eq:Alpha-variation-energy70}) and (\ref{eq:Alpha-variation-energy70-damping})
 the result is obtained.


\end{document}